\documentclass[11pt,a4paper]{article}

\usepackage{fullpage}

\usepackage{amssymb}
\usepackage{amsthm}
\usepackage{amsmath}            
\usepackage{latexsym}
\usepackage{epsfig}
\usepackage{amscd}
\usepackage{graphicx}
\usepackage[dvips]{color} 
\newcommand{\Real}{\mbox{$\mathbb{R}$}}         
\usepackage{hyperref}

\renewcommand{\proof}[1]{\noindent \normalfont{\textbf{Proof.} \  \  #1 \hfill $\Box$} \par}

\newcommand{\ben}{\begin{equation}}     
\newcommand{\eeqn}{\end{equation}}
\newcommand{\bey}{\begin{eqnarray}}
\newcommand{\eey}{\end{eqnarray}}

\newtheorem{thm}{Theorem}[section]
\newtheorem{prop}[thm]{Proposition}
\newtheorem{lemma}[thm]{Lemma}
\newtheorem{corol}[thm]{Corollary}
\newtheorem{defn}[thm]{Definition}

\begin{document}

\bibliographystyle{plain}

\title{The Connected Isentropes Conjecture in a Space of Quartic Polynomials}
\maketitle

Anca Radulescu, Mathematics Department, Stony Brook University

\begin{abstract}
\small{ This note is a shortened version of my dissertation thesis, defended at Stony Brook University in December 2004. It illustrates how dynamic complexity of a system evolves under deformations. The objects I considered  are quartic polynomial maps of the interval that are compositions of two logistic maps. In the parameter space $P^{Q}$ of such maps, I considered the algebraic curves corresponding to the parameters for which critical orbits are periodic, and I called such curves left and right bones. Using quasiconformal surgery methods and rigidity I showed that the bones are simple smooth arcs that join two boundary points. I also analyzed in detail, using kneading theory, how the combinatorics of the maps evolves along the bones. 
   The behavior of the topological entropy function of the polynomials in my family is closely related to the structure of the bone-skeleton. The main conclusion of the paper is that the entropy level-sets in the parameter space that was studied are connected.    }
\end{abstract}

\thispagestyle{empty} 
\def\IMSmarkvadjust{0 pt}
\def\IMSmarkhadjust{0 pt}
\def\IMSmarkhpadding{0 pt}
\def\IMSpubltext{Published in modified form:}
\def\SBIMSMark#1#2#3{
 \font\SBF=cmss10 at 10 true pt
 \font\SBI=cmssi10 at 10 true pt
 \setbox0=\hbox{\SBF \hbox to \IMSmarkhpadding{\relax}
                Stony Brook IMS Preprint \##1}
 \setbox2=\hbox to \wd0{\hfil \SBI #2}
 \setbox4=\hbox to \wd0{\hfil \SBI #3}
 \setbox6=\hbox to \wd0{\hss
             \vbox{\hsize=\wd0 \parskip=0pt \baselineskip=10 true pt
                   \copy0 \break%
                   \copy2 \break%
                   \copy4 \break}}
 \dimen0=\ht6   \advance\dimen0 by \vsize \advance\dimen0 by 8 true pt
                \advance\dimen0 by -\pagetotal
	        \advance\dimen0 by \IMSmarkvadjust
 \dimen2=\hsize \advance\dimen2 by .25 true in
	        \advance\dimen2 by \IMSmarkhadjust

%
%
  \openin2=publishd.tex
  \ifeof2\setbox0=\hbox to 0pt{}
  \else 
     \setbox0=\hbox to 3.1 true in{
                \vbox to \ht6{\hsize=3 true in \parskip=0pt  \noindent  
                {\SBI \IMSpubltext}\hfil\break
                \input publishd.tex 
                \vfill}}
  \fi
  \closein2
  \ht0=0pt \dp0=0pt
 \ht6=0pt \dp6=0pt
 \setbox8=\vbox to \dimen0{\vfill \hbox to \dimen2{\copy0 \hss \copy6}}
 \ht8=0pt \dp8=0pt \wd8=0pt
 \copy8
 \message{*** Stony Brook IMS Preprint #1, #2. #3 ***}
}

\SBIMSMark{2005/04}{May 2005}{}

\clearpage
\pagenumbering{arabic}
\setcounter{page}{0}
\tableofcontents

\clearpage

\section{Previous work and summary of results}

\indent

   This paper illustrates how dynamic complexity of a system evolves under deformations. This evolution is in general only partly understood. Attempts to give a quantitative approach have considered simple examples of dynamical systems and have made use of the topological entropy $h(f)$ as a particularly useful measure of the complexity of the iterated map $f$. There has been a lot of work about entropy; although not much on monotonicity.

   The logistic family $\{ f_{\mu}(x) = \mu x(1-x) \:,\: \mu \in [0,4] \}$
illustrates many of the important phenomena that occur in Dynamics. The theory in this case is the most complete (see [D]): $ \mu \rightarrow h(f_{\mu}) $ is continuous, monotonely increasing, and different values $h_{0} = h(f_{\mu})$ are realized for a single $\mu$ in some cases, but also for infinitely many in other cases. The cubic polynomials on the unit interval are organized as a 2-parameter family. In the compact parameter space of this family, the level sets of the entropy, called isentropes, were proved to be connected ([DGMT] and [MT]).

   In general, families of degree $d$ polynomials depend on $d-1$ parameters, so the same concepts are harder to inspect for higher degrees. It is most natural to research next a family of quartic polynomials that depends only on two parameters. This paper focuses on showing the \textbf{Connected Isentropes Conjecture} for the parameter space $P^{Q}$ of the family of alternate compositions of two logistic maps ([R]). The work is organized as follows:\\

  I briefly study the more general combinatorics of $2n$-periodic orbits under alternate iterations of two $(+,-)$ unimodal interval maps.

  I introduce a way to keep track of the succession of the orbit points along the unit interval $I$ by defining the \emph{order-data} as a pair of permutation $(\sigma,\tau) \in S_{n}^{2}$. If under alternate iterations of the two maps $h_{1}$ and $h_{2}$ the two critical orbits are periodic, their order-data turns out to be strongly connected to the kneading-data of the composition $h_{2} \circ h_{1}$.

   For a given order-data $(\sigma,\tau)$, I define the left/right bones in the parameter space $P^{Q}$ to be the subsets for which either critical point has periodic orbit of order-data $(\sigma,\tau)$. The bones are algebraic curves , and by definition left bones can only intersect right bones. A crossing is called a primary intersection if it corresponds to a pair of maps with common periodic bicritical orbit, and secondary intersection if it corresponds to a pair of maps with disjoint critical orbits.

   To obtain combinatorial properties of the bones, I compare the space $P^{Q}$ with a model space of compositions of stunted tent maps. This technique is not accidental; the stunted sawtooth maps are generally useful models in kneading-theory, because they are rich enough to encode in a canonical way all possible kneading-data of $m$-modal maps. The combinatorial results make crucial use of Thurston's Uniqueness Theorem, and of an extension of it due to Poirier , interpreted by [MT].\\

   In two following sections, I complete the description of the bones with two essential properties. 

   The bone-curves are ${\cal{C}}^{1}$-smooth and intersect transversally. Smoothness follows as in [M] at parameter points inside the hyperbolic components of $P^{Q}$. If the parameter point is outside these components, a quasiconformal surgery construction is necessary in order to perturb a map with a superattracting cycle to a map having an attracting  cycle with small nonzero  multiplier.

   The bones are simple arcs in $P^{Q}$ with two boundary points on $\partial P^{Q}$, in other words they contain no loops. [MT] proved the similar assertion in the case of cubic polynomials, either assuming true the well-known Fatou Conjecture or using a weaker theorem due to Heckman. I use instead a quite new and interesting rigidity result of [KSvS], that delivers density of hyperbolicity in my parameter space.

   I define the $n$-skeleton  $S_{n}^{Q}$ in $P^{Q}$ to be the union of all bones of period at most $2n$ , together with the boundary of the space. I put a dimension 2 topological cell structure on $P^{Q}$ as follows: the 0-cells are all intersections of bones in $S_{n}^{Q}$ and all boundary points of bones in $S_{n}^{Q}$;  the 1-cells are the 1-dimensional connected components obtained by deleting the 0-cells from the $n$-skeleton; the 2-cells are the 2-dimensional connected components of the complement of $S_{n}^{Q}$.\\

   The relations between entropy and the sequence of cell complexes is emphasized in the last section of the paper. If two points in $P^{Q}$ correspond to distinct values of the entropy, then any path connecting them crosses infinitely many bones. In more technical phrasing: for any $\epsilon >0$, there is a large enough $n$ for which the corresponding cell complex is fine enough to have variation of entropy less than $\epsilon$ on each of its closed cells. These considerations permit me to transport some topological properties of the isentropes from the previously mentioned model space to similar properties of isentropes in $P^{Q}$. More precisely, contractibility of isentropes in the stunted tent maps model space translates as connectedness of isentropes in $P^{Q}$.

\section{Combinatorics}  
\subsection{A discussion on the kneading-data} \label{section:subsection1} 
\indent
    
   Let $h:I \rightarrow I$ be an \emph{m-modal map} of the interval, i.e. there exist $\; 0<{\bf c_{1}} \leq {\bf c_{2}} \leq ... \leq {\bf c_{m}}<1 \;$  ``folding'' or ``critical points'' of $h$ such that $h$ is alternately increasing and decreasing on the intervals  $H_{0},...,H_{m}$ between the folding points.

   $$ I= \bigcup_{k=0}^{m} H_{k} \cup \bigcup_{j=1}^{m} \{ {\bf c_{j}} \} $$
   We say that $h$ is of shape $s=(+,-,+,...)$ if $h$ is increasing on $H_{0}$ and of shape $s=(-,+,-,...)$ if $h$ is decreasing on $H_{0}$. We say that $h$ is strictly $m$-modal if there is no smaller $m$ with the properties above.

   We define the itinerary $\Im(x)=(A_{0}(x),A_{1}(x),...)$ of a point $x \in I$ under $h$ as a sequence of symbols in $ {\cal{A}}=\{ H_{0},...,H_{m} \} \cup \{{\bf c_{1},...,c_{m}} \}$, where

   \[  \left\{ \begin{array}{ccc}
          \mbox{$A_{k}(x)=H_{j}$} & \mbox{, if} & \mbox{$f^{\circ k}(x) \in H_{j}$} \\
          \mbox{$A_{k}(x)={\bf c_{j}}$} & \mbox{, if} & \mbox{$f^{\circ k}(x)={\bf c_{j}}$}
          \end{array} 
       \right. \]

   The kneading sequences of the map $h$ are defined as the itineraries of its folding values:

   $$ {\cal{K}}_{j}={\cal{K}}({\bf c_{j}})=\Im(f({\bf c_{j}})), \; j=\overline{1,m-1}$$

   The kneading-data {\bf{K}} of $h$ is the $m$-tuple of kneading-sequences:

   $${\bf K}=({\cal{K}}_{1},...,{\cal{K}}_{m})$$

   The simplest example of an $m$-modal map is a \emph{sawtooth map} with $m$ teeth (see figure 1.1(a)).

\begin{figure}[h]
\begin{center}
\input{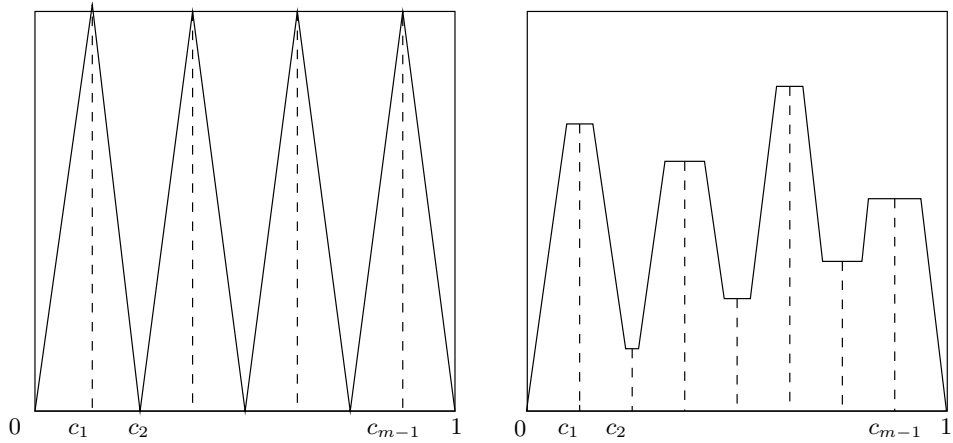}\\
\caption{\emph{(a)Sawtooth map of the interval. (b)Stunted sawtooth map}}
\end{center}
\end{figure}

 We call a \emph{stunted sawtooth map} a sawtooth map whose vertexes have been stunted by plateaus placed at chosen heights (see figure 1.1(b)).Its critical points are considered to be the centers of the plateaus. In the next sections we will focus our attention specifically on tent maps ($1$-modal sawtooth maps) and on their stunted version, which we will call stunted tent maps.\\

   Another simple and rich example of $m$-modal maps is the collection of $(m-1)$-degree polynomials from $I$ to itself. The ``folding points'' could be taken in this case to be the critical points of the polynomial (in the classical sense) of odd order. In the context of polynomial $m$-modal maps, we have a powerful tool to use in the statement of \emph{Thurston's Uniqueness Theorem}.

   \begin{defn}  \emph{A polynomial map is called post-critically finite if the orbit of every critical point is periodic or eventually periodic.}
   \end{defn}

   \begin{thm}\textbf{Thurston Uniqueness Theorem for Real Polynomial Maps:}       A post-critically finite real polynomial map of degree m+1 with m distinct real critical points is uniquely determined, up to a positive affine conjugation, by its kneading data.
   \end{thm}

   We will also use a converse of this basic theorem of Thurston, due to Poirier (as interpreted by ~\cite{MT}).\\

   \begin{defn} \emph{ We say that a symbol sequence $\Im(x)=(A_{0}(x),A_{1}(x),...)$ is} flabby \emph{if some point of the associated orbit which is not a folding point has the same itinerary as an immediately adjacent folding point. A symbol sequence is called} tight \emph{if it is not flabby. The kneading data of a map is tight if each of its kneading sequences is tight.}
    \end{defn}

   \begin{lemma} The kneading data of a stunted sawtooth map is tight if and only if the orbit of each folding point never hits a plateau except at its critical point.
   \end{lemma}

   \begin{thm} Suppose that the m-modal kneading data ${\bf K}$ is admissible for some shape $s$, with ${\cal}{K}_{i} \neq {\cal}{K}_{j}$ for all $i$. There exists a post-critically finite polynomial map of degree m+1 and shape $s$ with kneading-data ${\bf K}$ if and only if each ${\cal}{K}_{i}$ is periodic or eventually periodic, and also tight. This polynomial is always unique when it exists, up to a positive affine change of coordinates, or as a boundary anchored map of the interval.
   \end{thm}

\subsection{Definitions and first goals}\label{section:subsection2}
\indent

   In the light of the general definition given in section 2.1, a \emph{ boundary anchored, (+,-) unimodal map} of the unit interval is a $h:I=[0,1] \rightarrow I$ such that $h(0)=h(1)=0$ and such that there exists $\gamma \in (0,1)$, called \emph{folding} or \emph{critical} point, with $h$ increasing on $(0,\gamma)$ and decreasing on $(\gamma,1)$. The orbit of a point $x \in I$ under a such $h$ will be the sequence of iterates $(h^{\circ n}(x))_{n \geq 0}$. The itinerary of $x$ under $h$ is the sequence $(J_{0},J_{1},...)$ of symbols $L$ (left), $R$ (right) and $\Gamma$ (center, or critical) such that:

       \[ \left\{ \begin{array}{ccc}
       \mbox{$J_{j} = L ,$} & \mbox{if $h^{\circ j}(x) < \gamma$ }\\
       \mbox{$J_{j} = R ,$} & \mbox{if $h^{\circ j}(x) > \gamma$ }\\
       \mbox{$J_{j} = \Gamma ,$} & \mbox{if $h^{\circ j}(x) = \gamma$ }
               \end{array}
       \right. \]

   The next few sections of this paper are dedicated to study the combinatorics of the dynamical system I am considering: generated by alternate iterates of two unimodal interval maps. In this sense, it is convenient to consider two copies of the unit interval $I_{1}=I_{2}=I$ and think of our pair of maps $(h_{1},h_{2})$ as a self map of the disjoint union $I_{1} \sqcup I_{2} \rightarrow I_{1} \sqcup I_{2}$, which carries $I_{1}$ to $I_{2}$ as $h_{1}$ and $I_{2}$ to $I_{1}$ as $h_{2}$, with critical points $\gamma_{1} \in I_{1}$ and $\gamma_{2} \in I_{2}$, respectively. \\
 
   We call an \emph{orbit} under the pair $(h_{1},h_{2})$ a sequence: 

   $$x \rightarrow h_{1}(x) \rightarrow h_{2}(f_{1}(x)) \rightarrow h_{1}(h_{2}(h_{1}(x)))...$$
We say a such orbit is \emph{critical} if it contains either critical point $\gamma_{1}$ or $\gamma_{2}$ and we say it is \emph{bicritical} if it contains both. We call the \emph{itinerary} of a point $x$ under $(h_{1},h_{2})$ the infinite sequence $\Im(x)=(J_{k}(x))_{k \geq 0}$ of alternating symbols in $\{ L_{1},\Gamma_{1},R_{1} \}$ and $\{ L_{2},\Gamma_{2},R_{2} \}$ that expresses the positions of the iterates of $x$ in $I_{1}$ and $I_{2}$ with respect to $\gamma_{1}$ or $\gamma_{2}$.

   Clearly, not all arbitrary symbol sequences are in general admissible as itineraries of a point under a pair of given maps. 

   It is fairly easy to show that for a fixed pair $(h_{1},h_{2})$ of $(+,-)$ unimodal maps, the pair of critical itineraries $(\Im(\gamma_{1}),\Im(\gamma_{2}))$ determines the kneading-data of $h_{2} \circ h_{1}$ and conversely. In particular this applies to pairs of stunted tent maps and to pairs of logistic maps (which are the object of this paper).

   For a given pair of maps $(h_{1},h_{2})$, I will use the regular total order on admissible itineraries (see ~\cite{CE}), which is consistent with the order of points on the real line:

   $$ \Im(x) < \Im(x') \Rightarrow x < x' $$
   $$ x < x' \Rightarrow \Im(x) \leq \Im(x')$$

    We say that the orbit of $x$ is periodic of period $2n$ under $(h_{1},h_{2})$ if $n$ is the smallest positive integer such that $(h_{2} \circ h_{1})^{\circ n}(x)=x$ (i.e. $x$ has period $n$ under the composition $(h_{2} \circ h_{1})$) . I will use the following notation for a $2n$-periodic orbit under $(h_{1},h_{2})$:

\begin{equation} x_{1}=x_{i_{1}} \stackrel{h_{1}}{\longrightarrow} y_{j_{1}} \stackrel{h_{2}}{\longrightarrow} x_{i_{2}} \stackrel{h_{1}}{\longrightarrow} ... \stackrel{h_{1}}{\longrightarrow} y_{j_{n}} \stackrel{h_{2}}{\longrightarrow} x_{i_{1}}
\end{equation}

\noindent
where $(x_{i})_{i=\overline{1,n}} \subset I_{1}$ and $(y_{j})_{j=\overline{1,n}} \subset I_{2}$ are both increasing.

  \begin{defn} \emph{The \emph{order-data} of the periodic orbit (1) is the pair $(\sigma,\tau)$ of permutations in $S_{n}$ given by:} 

    $$ h_{1}(x_{i})=y_{\sigma_{i}} \: \,$$ 
    $$ h_{2}(y_{j})=x_{\tau_{j}} \: , $$
\emph{so that $\sigma_{i_{k}}=j_{k}$ and $\tau_{j_{k}}=i_{k+1}$. (Here the subscripts must be understood as integers mod $n$, e.g. $i_{n+1}=i_{1}=1$.)}

  \end{defn}

   An admissible order-data is a $(\sigma,\tau) \in S_{n}^{2}$ which is achieved as order-data of a periodic orbit of some pair $(h_{1},h_{2})$ of interval unimodal maps.\\

   The (+,-) unimodal shape of $h_{1}$ and $h_{2}$ imposes a set of necessary and sufficient conditions for a ($\sigma, \tau$) to be ``admissible'':

    \[ (I) \left  \{ \begin{array}{ll}
         \mbox{If $\sigma_{i+1} < \sigma_{i}$ ,} & \mbox{then $\sigma_{j+1} < \sigma_{j} ,  \forall j \geq i $} \\
         \mbox{If $\tau_{i+1} < \tau_{i}$ ,} & \mbox{then $\tau_{j+1} < \tau_{j} , \forall j \geq i $} 
         \end{array} \right. \] 
   
   $ \hspace{20mm} (II) \quad \tau \circ \sigma$ is a cyclic 
   permutation (i.e. has no smaller cycles).\\

\begin{figure}[h]
\begin{center}
\input{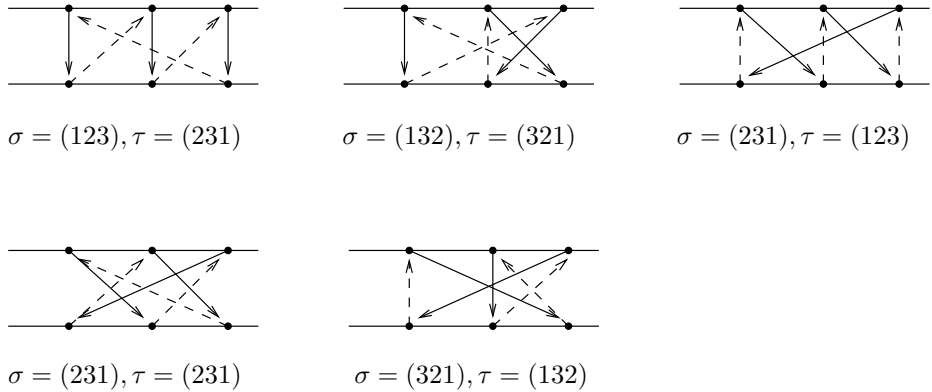}\\
\caption{\emph{All admissible order-data $(\sigma,\tau)$of period $2n=6$. Each schetch represents the interval $I=I_{1}$ on top, with the orbit points $x_{1}<x_{2}<x_{3}$ and the interval $I=I_{2}$ undernieth, with the orbit points $y_{1}<y_{2}<y_{3}$.}}
\end{center}
\end{figure}

   A first goal will be to research the relation between the itinerary and the order-data of a periodic critical orbit.

    Suppose $\gamma_{1}$ is periodic  of period $2n$ under $(h_{1},h_{2})$ and let $\quad  x_{1}=x_{i_{1}} \rightarrow y_{j_{1}} \rightarrow ... \rightarrow y_{j_{n}} \rightarrow x_{i_{1}} \quad $ be its orbit. Then the order-data $(\sigma,\tau) \in S_{2n}^{2}$ of the orbit determines its itinerary via the position of the element $y_{j_{l}} \in I_{2}$ closest to $\gamma_{2}$. In other words, there are at most two critical itineraries corresponding to a given order-data.  If in particular the orbit is bicritical, then $y_{j_{l}}=\gamma_{2}$ and the itinerary is completely defined. 

    Note also that the order of points in a critical periodic orbit of a $(+,-)$ unimodal map is \emph{strictly} preserved in the order of their itineraries, i.e. $x < x'$ implies $\Im(x) < \Im(x')$. Hence conversely, knowing the itinerary $\Im$ of the bicritical orbit, we can obtain the order of occurence of the orbit points in $I_{1}$ and $I_{2}$. This proves the following:

   \begin{thm} If the orbit of $\gamma_{1}$ is bicritical of period 2n under a pair of (+,-) unimodal maps $(h_{1},h_{2})$, then the itinerary of $\gamma_{1}$ determines the order-data of the orbit and conversely.
   \end{thm}

\subsection{Parameter spaces} \label{section:subsection3}
\indent

   We plan to study in more detail the dynamics of a particular family of such pairs of interval unimodal maps that we have generally described in the previous section.

   Recall that the logistic map (with critical value $v$) is defined as $q_{v}(x)=4vx(1-x), \: x \in \mathbb{R}$. Clearly $q_{v}(0)=q_{v}(1)=0$, for any value of the parameter $v$. Moreover, for values of $v \in [0,1]$, $q_{v}$ carries the unit interval to itself, so it is a boundary anchored, $(+,-)$ unimodal interval map. 

   Our goal is to study the dynamics of compositions of pairs of such maps: $q_{w} \circ q_{v}$, where $(v,w) \in [0,1]^{2}$. We call the family of pairs $(q_{v},q_{w})$ of logistic maps of the unit interval the \textbf{$Q$-family}, and we parametrize it by the pair of critical values, so that the parameter space will be: 

  $$P^{Q}=\{ (v,w) \in [0,1] \times [0,1] \}=[0,1]^{2}$$

   The behavior of the pairs in the $Q$-family is not very well-understood. We will compare it to the dynamics in a ``model'' family much easier to research, the family of pairs of stunted tent maps:

    $$ st_{v}:I_{1} \rightarrow I_{2} ,\; \gamma_{1}=\frac{1}{2} $$ 
    $$ st_{w}:I_{2} \rightarrow I_{1} ,\; \gamma_{2}=\frac{1}{2} $$
where 

     \[ st_{v}(x) =\left\{\begin{array}{lll} 
                  2x & \mbox{if $x \leq \frac{v}{2}$} \\
                   v & \mbox{if $\frac{v}{2} \leq x \leq 1-\frac{v}{2}$} \\
                  2-2x & \mbox{if $x \geq 1-\frac{v}{2}$}
                 \end{array}
                 \right. \]

   Recall that the ``critical point'' of such a stunted tent map was taken by convention to be the midpoint $\gamma=\frac{1}{2}$, hence the critical value is $st(\gamma)=v$. We call the family of pairs of such maps the ST-family. Its corresponding parameter space will be denoted by:

   $$P^{ST}= \{ (v,w) \in [0,1] \times [0,1] \}$$

\vspace{3mm}

   We aim to obtain dynamical results in $P^{Q}=[0,1]^{2}$. However, proving similar results in the parameter space $P^{ST}=[0,1]^{2}$ of ``approximating'' stunted tent maps would be a good start. Comparison of the two spaces will be a strategy very frequently used within the combinatorics sections. The strong topological correspondence between the two families will eventually be sustained with a rigurous proof and will enable us to translate topological properties from one to the other. 
\subsection{The combinatorics in the ST-family}\label{section:subsection4}
\indent

    I will focus next on how the combinatorics in section 2.2 applies to the model family I am interested in, namely the $ST$-family:

    \begin{thm} Given $(\sigma,\tau) \in S_{n}^{2}$ admissible order-data, there is a unique pair of stunted tent maps $(st_{v},st_{w})$ with periodic bicritical orbit of order-data $(\sigma,\tau)$.
   \end{thm}

   \proof{ Let $\Im$ be a sequence of alternating symbols in $\{L_{1},R_{1},\Gamma_{1} \}$ and $\{ L_{2},R_{2},\Gamma_{2} \}$, admissible as a bicritical itinerary of period $2n$ under a pair of unimodal maps:

     $$ \Im = (J_{0}=\Gamma_{1},J_{1},J_{2},...,J_{2l},J_{2l+1}=\Gamma_{2},J_{2l+2},...,J_{2n-1},J_{2n}=\Gamma_{1},... ) $$
where $J_{2n+k}=J_{k}$ for all $k$ and $J_{k} \neq \Gamma_{1}, \Gamma_{2}$, for all $k$ nonequivalent to $1,...,2l$ mod $2n$. There exists a unique pair of stunted tent maps ($st_{v},st_{w}$) that has a bicritical orbit of period $2n$:
 
  $$ x_{1}=x_{i_{1}} \rightarrow y_{j_{1}} \rightarrow ... y_{j_{n}} \rightarrow 
     x_{i_{1}}=x_{1} $$ 
having $\Im$ as its itinerary.

    To prove the existence, it is easier to consider an orbit through a pair of tent maps that has the respective itinerary, then stunt the maps at the highest values of the orbit in $I_{1}$ and $I_{2}$, respectively. The uniqueness follows: starting with the critical points $\gamma_{1}$ and $\gamma_{2}$, iterate backwards using the itinerary $\Im$ to obtain the values of $v$ and $w$.

    Going back to the proof of our theorem: given  an admissible order-data ($\sigma, \tau$) $\in S_{n}^{2}$ for a required bicritical orbit, we can determine the itinerary $\Im$ of the orbit. As shown above, we can find a unique pair ($st_{v},st_{w}$) of stunted tent maps with a bicritical orbit of length $2n$ and itinerary $\Im$. By Theorem 1.3.2, the order-data for the orbit we have found will be ($\sigma,\tau$). }

\vspace{3mm}

     To make the discussion a step more general, we look next at pairs of arbitrary unimodal maps for which both critical points $\gamma_{1}$ and $\gamma_{2}$ are periodic. There are two possible cases that can occur: a bicritical orbit (discussed earlier) and two disjoint critical orbits.

   \begin{defn} \emph{Let $(\sigma,\tau) \in S_{m+n}^{2}$ be a pair of permutations decomposable into two cycles: $(\sigma_{1},\tau_{1}) \in S_{m}^{2}$ and $(\sigma_{2},\tau_{2}) \in S_{n}^{2}$. We say that two disjoint periodic orbits $o_{1}$ and $o_{2}$ under a pair $(h_{1},h_{2})$ of (+,-) unimodal maps have} joint order-data \emph{$(\sigma,\tau)$ if:}
  
   \begin{enumerate}

   \item   \emph{$o_{1}$ has order-data $(\sigma_{1},\tau_{1})$ and $o_{2}$ has order-data $(\sigma_{2},\tau_{2})$;}
   
   \item   \emph{ the order of the points in $I_{1}$ and $I_{2}$ is given by $(\tau \circ \sigma)$ and $(\sigma \circ \tau)$ respectively. (see ``order-type'' ~\cite{MT})}\\
   \end{enumerate}
   
   \end{defn}

     We will say about a permutation $(\sigma,\tau) \in S_{m+n}$ that it is ``admissible'' as a joint order-data, if there exist two disjoint orbits under some pair of (+,-) unimodal maps which have joint order-data $(\sigma,\tau)$.

    Similarly as for regular order-data, one can obtain the following two results:

   \begin{thm} Let $o_{1}$ and $o_{2}$ be disjoint critical orbits under a pair $(h_{1},h_{2})$ of (+,-) unimodal maps. Their itineraries determine their joint order-data and conversely.
   \end{thm}

    \begin{thm}  Given $(\sigma,\tau)=((\sigma_{1},\tau_{1}), (\sigma_{2},\tau_{2})) \in S_{m+n}^{2}$ admissible joint order-data, there exists a unique pair $(st_{v},st_{w})$ of stunted tent maps with disjoint critical orbits $o_{1} \ni \gamma_{1}$ and $o_{2} \ni \gamma_{2}$  having joint order-data $(\sigma, \tau)$.
   \end{thm}

\subsection{Description of bones in the ST-family} \label{section:subsection5}

\indent

  \begin{defn} \emph{Fix an admissible order-data $(\sigma,\tau) \in S_{n}^{2}$. By a \emph{left bone} in the parameter space $I_{2} \times I_{1}$ for the $ST$-family we mean the set of pairs $(v,w) \in I_{2} \times I_{1} = [0,1]^{2}$ such that the critical point $\gamma_{1} \in I_{1}$ has under $(st_{v},st_{w})$ a periodic orbit of given period $2n$ and given order-data $(\sigma,\tau)$.} 
   \end{defn}

   We will use the notation $B_{L}^{ST}(\sigma,\tau)$, or $B_{L}^{ST}$ if there is no ambiguity. We define a \emph{right bone} symmetrically (i.e. we require $\gamma_{2}$ to be periodic of specified period and order-data) and we denote it by $B_{R}^{ST}(\sigma,\tau)$, or $B_{R}^{ST}$. We will need later a more comprehensive approach to the left and right bones and their properties.

   Recall (from theorem 2.8) that: \emph{ There is a unique pair $(v_{0},w_{0}) \in B_{L}^{ST}$ such that the periodic orbit of $\gamma_{1}$ is bicritical (i.e. hits $\gamma_{2}$) under $(st_{v_{0}},st_{w_{0}})$.}\\

   \begin{thm} For each admissible order-data $(\sigma,\tau)$, let $(v_{0},w_{0})$ be the parameter pair for the associated bicritical orbit in the $ST$-family. Then there are unique numbers $v_{1}<v_{0}<v_{2}$ so that the left bone $B_{L}^{ST}(\sigma,\tau)$ is the union $\{ v_{1},v_{2} \} \times [w_{0},1] \cup (v_{1},v_{2}) \times \{ w_{0} \}$ of three line segments, as illustrated in figures 3 and 4. The description of the right bone $B_{R}^{ST}(\sigma,\tau)$ is completely analogous.
   \end{thm}

\begin{figure}[h]
\begin{center}
\input{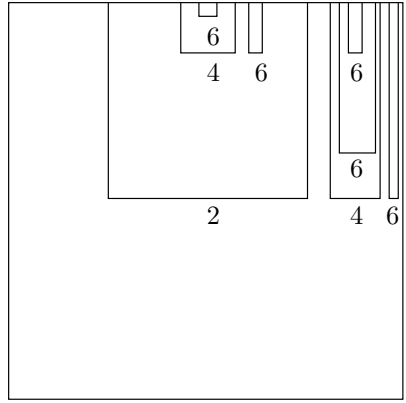}\\
\caption{\emph{Left bones in the ST-family of period at most 6.We marked by (2) the unique bone of period 2, corresponding to order-data in $(\sigma=(1),\tau=(1)) \in S_{1}^{2}$. (4) are the 2 bones of period 4 and having the two possible order-data $(\sigma=(12),\tau=(1)(2))$ or $(\sigma=(1)(2),\tau=(12)) \in S_{2}^{2}$. (6) are the bones of period 6 and one of the 5 admissible order-data: $(\sigma=(123),\tau=(231)$, $(\sigma=(132),\tau=(321)),(\sigma=(231),\tau=(231)), (\sigma=(321),\tau=(132))$ or $(\sigma=(231),\tau=(123))$ }}
\end{center}
\end{figure}

   We can determine the shape of $B_{L}^{ST}$, hence prove 2.13, by constructive means, starting with the point $(v_{0},w_{0})$.

   Under $(st_{v_{0}},st_{w_{0}})$: $ \gamma_{1} \rightarrow v_{0} \rightarrow ... \rightarrow \gamma_{2} \rightarrow w_{0} \rightarrow ... \rightarrow \gamma_{1} $. The bicritical orbit only hits each plateau once, at its center.

   By sliding the first plateau up and down, the orbit of $\gamma_{1}$ will change in a continuous way. For a fixed height $v$ of the first plateau, call $y_{l}(v)$ the element in $I_{2}$ closest to $\gamma_{2}$ in the orbit of $\gamma_{1}$ under $(st_{v},st_{w})$. Clearly, if $v=v_{0}$, then $y_{l}(v_{0})=\gamma_{2}$.

   We can move $v$ continuously within an interval $[v_{1},v_{2}]=[v_{0}-\epsilon,v_{0}+\epsilon], \; \epsilon>0$ such that $y_{l}(v)$ moves from $w_{0}/2$ to $1-w_{0}/2$. Along the process, the orbit stays periodic and the order of the occurrence of points remains consistent with $(\sigma,\tau)$.

   It is not hard to see that $B_{L}^{ST}=\sqcup=\sqcup(\sigma,\tau)=\{ v_{1},v_{2} \} \times [w_{0},1] \cup (v_{1},v_{2}) \times \{ w_{0} \}$.

   In particular, there are exactly two values $v=v_{1}$ and $v=v_{2}$ such that the orbit of $\gamma_{1}$ has given order-data $(\sigma,\tau)$ under $(st_{v},st_{1})$ (i.e. there are exactly two points of $B_{L}^{ST}$ on $[0,1] \times \{ 1 \}$).

\begin{figure}[h]
\begin{center}
\input{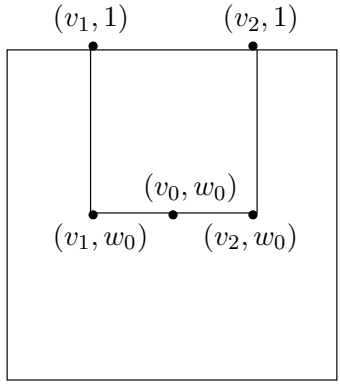}\\
\caption{\emph{ $B_{L}^{ST}=\sqcup=\{ v_{1},v_{2} \} \times [w_{0},1] \cup (v_{1},v_{2}) \times \{ w_{0} \} \ni (v_{0},w_{0})$}}
\end{center}
\end{figure}

\subsection{Important points on the bones} \label{section:subsection6}
\indent

   We aim to compare the parameter spaces for the two families: the $Q$-family and the $ST$-family.
   
   In either space, we consider the left and right $2n$-bones for a given admissible $(\sigma,\tau) \in S_{n}^{2}$:\\

   $B_{L}(\sigma,\tau)$= the set of all parameters for which $\gamma_{1}$ has periodic orbit of order-data $(\sigma,\tau)$ under the respective pair of maps \\

   $B_{R}(\sigma,\tau)$= the set of all parameters for which $\gamma_{2}$ has periodic orbit of order-data $(\sigma,\tau)$ under the pair of maps\\

    For any fixed admissible $(\sigma,\tau) \in S_{n}^{2}$, I will call the bones in $P^{ST}$: $B^{ST}_{L}$, $B^{ST}_{R}$ and the ones in $P^{Q}$: $B^{Q}_{L}$, $B^{Q}_{R}$.\\

    \textbf{Remarks}: (1) In either parameter space, any two left bones are disjoint and any two right bones are disjoint by definition.
 
    (2) It follows easily from theorem 2.13 that two bones in the ST-family can cross only in 0,2 or 4 points.\\

   \begin{defn} \emph{In either parameter space, an intersection of $B_{L}(\sigma_{1},\tau_{1})$ and $B_{R}(\sigma_{2},\tau_{2})$ is called a \emph{primary intersection} if $(\sigma_{1},\tau_{1})=(\sigma_{2},\tau_{2})$ and there is a bicritical orbit with this order-data under the pair of maps. It is called a \emph{secondary intersection} if the two critical orbits are disjoint, of distinct order-data $(\sigma_{1},\tau_{1})$ and respectively $(\sigma_{2},\tau_{2})$, and joint order-data $(\sigma,\tau)$.}
    \emph{ A \emph{capture point} on $B_{L}(\sigma_{1},\tau_{1})$ in either $P^{ST}$ or $P^{Q}$ is a pair of maps for which $\gamma_{2}$ eventually maps on $\gamma_{1}$ such that it has an eventually periodic, but not periodic, orbit. We define symmetrically a capture point on $B_{R}(\sigma_{2},\tau_{2})$.}

   \end{defn}

\vspace{3mm}

\begin{figure}[h]
\begin{center}
\input{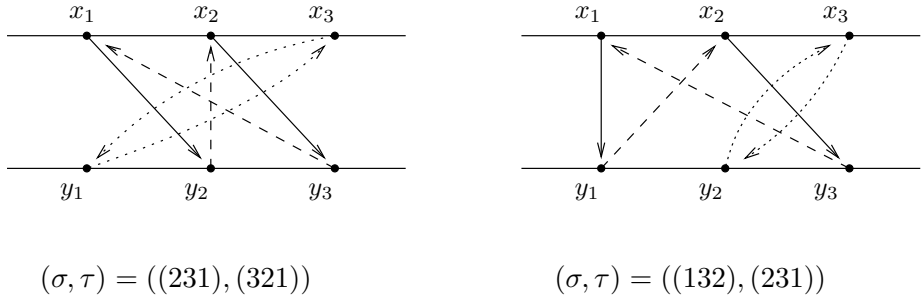}\\
\caption{\emph{Combinatorics of the two secondary intersections of a period 4 left bone with a period 2 right bone.}}
\end{center}
\end{figure}

\begin{figure}[h]
\begin{center}
\input{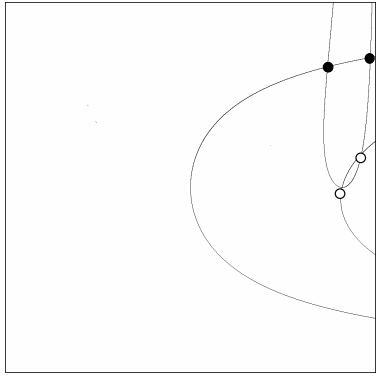}\\
\caption{\emph{ The left $4$-bone of order-data $(\sigma,\tau)=((1,2),(2,1))$ crosses the right $2$-bone at two secondary intersections with joint order-data $((231),(321))$ and $((132),(231))$ (filled dots, also see figure 5) and crosses the corresponding right 4-bone at a primary intersection with order-data $(\sigma,\tau)$ and at a secondary intersection with joint order-data $((1243),(3421))$ (empty dots).}}
\end{center}
\end{figure}

   Theorem 2.8 equipped us with a bijection between admissible order-data and primary intersections in $P^{ST}$. Theorem 2.11 extended the result with a bijection between admissible joint order-data and secondary intersections. The next statement is a further extension for capture itineraries and can be proved similarly with the direct implication in theorem 2.11.

   \begin{thm} Suppose the two critical points of a pair of unimodal maps are such that one of them has a closed orbit and the other maps on this closed orbit after a finite number of iterates, but without being periodic itself. Let $\Im_{1}$ and $\Im_{2}$ be the itineraries of the two critical points. Then there exists at least a pair $(st_{v},st_{w})$ with critical itineraries $\Im_{1}$ and $\Im_{2}$, respectively. (i.e.: There exists at least a capture point in $P^{ST}$ with given ``capture'' critical itineraries.) 
 
   \end{thm}

 \subsection{More on kneading-data} \label{section:subsection7} 
\indent

   In this section we will construct a bijective correspondence of bones intersections between our two parameter spaces $P^{ST}$ and $P^{Q}$. For the proof, it is necessary to view the composition $q_{w} \circ q_{v}$ of two logistic maps either as a 3-modal map with three critical points in $I=I_{1}$: ${\bf c_{1}} \leq {\bf c_{2}} \leq {\bf c_{3}}$, with ${\bf c_{2}}=\gamma_{1}$ and $q_{v}({\bf c_{1}})=q_{v}({\bf c_{3}})=\gamma_{2}$ or as a unimodal map with folding point $\gamma_{1}$, in case $q_{v}(x)=\gamma_{2}$ has a double real root or two complex roots. I will use rigidity theorems that involve essentially properties of the kneading-data. \\

   Let us look in more detail at the possible kneading-data of the maps in $P^{ST}$ and $P^{Q}$.\\

   \textbf{Maps in $P^{ST}$}: For any $(v,w)\in P^{ST}$, the map $st_{w} \circ st_{v}$ could be considered 3-modal, with folding points $c_{1}=\frac{1}{4}, c_{2}=\gamma_{1}=\frac{1}{2}$ and $c_{3}=\frac{3}{4}$.

   $${\cal{A}}^{ST}=\{ [0,\frac{1}{4}),\frac{1}{4},(\frac{1}{4},\frac{1}{2}),\frac{1}{2},(\frac{1}{2},\frac{3}{4}),\frac{3}{4},(\frac{3}{4},1] \}$$.

   and

   $${\bf K_{ST}}=({\cal{K}}(c_{1}),{\cal{K}}(c_{2}),{\cal{K}}(c_{3}))$$

   We can consider $P^{ST}$ as made of three parts: $P^{ST}=P^{ST}_{1} \cup P^{ST}_{2} \cup P^{ST}_{3}$, where $P_{1}^{ST}=\{ (v,w) \in [0,1]^{2}, \; w \geq 2v \}$, $P_{2}^{ST}=\{ (v,w), \; w<2v, \; w \leq 2-2v \}$ and $P_{3}^{ST}=\{ (v,w),\; w>2-2v \}$. \\

   \textbf{I.} Clearly there are no right bones in $P^{ST}_{1}$, hence no bones intersections.

   \textbf{II.} $P^{ST}_{2}$ contains no secondary intersections, since $\frac{w}{2} \leq v \leq 1-\frac{w}{2}$, so $ st_{w}(\gamma_{2})=w=(st_{w} \circ st_{v})(\gamma_{1})$.

   Moreover, if $(v,w) \in P^{ST}_{2}$ is a primary intersection, then the map $st_{v} \circ st_{w}$ is strictly 3-modal, with only one exception: $(v,w)=(\frac{1}{2},\frac{1}{2})$.

    \textbf{III.} For $(v,w) \in P^{ST}_{3}$ we clearly have that $st_{w} \circ st_{v}$ is strictly 3-modal, hence ${\cal{K}}(c_{1})={\cal{K}}(c_{3}) \neq {\cal{K}}(c_{2})$.\\

   \textbf{Maps in $P^{Q}$}: The behavior of the degree 4 polynomials in the $Q$-family is also different for distinct values of the parameters.\\

\begin{figure}[h]
\begin{center}
\input{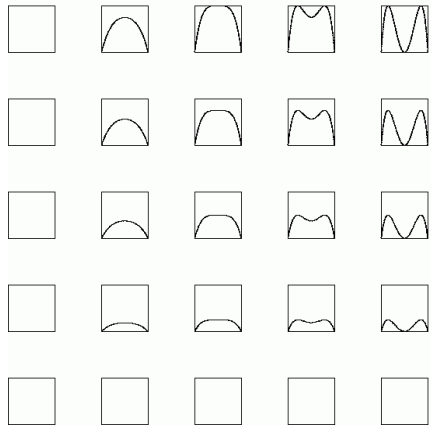}\\
\caption{\emph{ A few examples of behavior of maps in $P^{Q}$. The critical points of the quartic map $q_{w} \circ q_{v}$ are distinct and real for $v >\frac{1}{2}$,all coincide for $v=\frac{1}{2}$, while two of them are complex for $v<\frac{1}{2}$.}}
\end{center}
\end{figure}

   \textbf{I.} If $v<\frac{1}{2}$, then $q_{w} \circ q_{v}$ has only one real critical point $C_{2}=\gamma_{1}=\frac{1}{2}$ and two complex $C_{1},C_{3} \in \mathbb{C} \backslash \mathbb{R}$. 
   
   This parameter subset will be of somewhat less interest, as it is not crossed by any right bones, hence contains no bones intersections. Indeed, if $q_{v}(x) < \frac{1}{2}, \; \forall x \in I_{1}$, no orbit can go through $\gamma_{2}$.\\

   \textbf{II.} If $v=\frac{1}{2}$, then $q_{w} \circ q_{v}$ has a degenerate real critical point $C_{1}=C_{2}=C_{3}=\gamma_{1}$. This line contains primary intersections with right bones. More precisely, if a left bone hits $\{ v=\frac{1}{2} \}$, then the crossing point is its primary intersection. However, in this case $q_{v} \circ q_{w}$ is strictly 3-modal, with the exception of $v=w=\frac{1}{2}$, which is the period 2 primary intersection.\\

   \textbf{III.} If $v>\frac{1}{2}$, there are three distinct real critical points for $q_{w} \circ q_{v}$:

   $C_{1}<C_{2}=\gamma_{1}<C_{3}$, with $q_{v}(C_{1})=q_{v}(C_{3})=\gamma_{2}$

   The map is 3-modal:

    $${\cal{A}}^{Q}=\{ [0,C_{1}),C_{1},(C_{1},C_{2}),C_{2},(C_{2},C_{3}),C_{3},(C_{3},1] \}$$.

   and

   $${\bf K_{Q}}=({\cal{K}}(C_{1}),{\cal{K}}(C_{2}),{\cal{K}}(C_{3}))$$ \\

\noindent \textbf{Remark.} We emphasize that $q_{w} \circ q_{v}$ has complex critical points iff $v<\frac{1}{2}$. If the point $(v,w)$ is on a bone, it cannot be in the region $\{ v<\frac{1}{2}, \, w<\frac{1}{2} \}$, so $w \geq \frac{1}{2}$. Hence in this case the map $q_{v} \circ q_{w}$ corresponding to the symmetric point $(w,v)$ on the corresponding right bone has real critical points, non-degenerate if $w \neq \frac{1}{2}$.\\

   A correspondence is already apparent between the shape and position of two left bones with identical order-data in the two spaces $P^{ST}$ and $P^{Q}$. For instance, the unique primary intersection of period two: $(v,w)=(\frac{1}{2},\frac{1}{2}) \in P^{ST}$ clearly corresponds combinatorialy to the identical point $(v,w)=(\frac{1}{2},\frac{1}{2}) \in P^{Q}$. We will consider at least this case classified in our future analysis. The following theorems will therefore concern specifically the strictly 3-modal case ( applicable for either $st_{w} \circ st_{v}$ and $q_{w} \circ q_{v}$ or $st_{v} \circ st_{w}$ and $q_{v} \circ q_{w}$).\\

\subsection{The correspondence of the  \\ bones intersections} \label{section:subsection8}
\indent

   I use Thurston's Theorem and its extension for boundary anchored polynomials of degree four and shape (+,-,+,-) to construct in this section a bijection between bones crossings in the two parameter spaces. For the rest of section 2, we will adapt our notation to distinguish between parameters $(v,w) \in P^{ST}$ and parameters $(v',w') \in P^{Q}$.

  \begin{thm}  Let $(\sigma,\tau) \in S_{n}^{2}$ be admissible order-data. There is a unique primary intersection $(v',w')$ in $P^{Q}$ with this data and conversely.
  \end{thm}

   \proof{ \textbf{Uniqueness}: Suppose we have a pair $(v,w) \in P^{Q}$ with a bicritical orbit of order-data $(\sigma,\tau)$. We implicitly know the itinerary of the bicritical orbit, hence the kneading sequences of the three real distinct critical points $C_{1}< C_{2}=\frac{1}{2}<C_{3}$ of $q_{w} \circ q_{v}$ ( if $v>\frac{1}{2}$) or $q_{v} \circ q_{w}$ (if $w>\frac{1}{2}$). By Thurston's Theorem, the boundary anchored polynomial of degree 4 with the expected kneading data is unique, implying the uniqueness of the pair $(q_{v},q_{w})$ with the given order-data.

   \textbf{Existence}: Let $(st_{v},st_{w})$ be the pair of stunted tent maps with bicritical orbit of order-data $(\sigma,\tau)$. We know by theorem 2.7 that we can determine the itinerary of this bicritical orbit. If we exclude the case $v=w=\frac{1}{2}$, which is already classified, then either $st_{w} \circ st_{v}$ or $st_{v} \circ st_{w}$ is strictly 3-modal (say $st_{w} \circ st_{v}$, to fix our ideas). We know the kneading-data \textbf{K} for $st_{w} \circ st_{v}$, which should also be the kneading-data for the polynomial $q_{w'} \circ q_{v'}$ that we want to find. We hence need to prove existence of a polynomial of degree 4 with the required kneading-data \textbf{K} and then show that it can be written as a composition of two logistic maps $q_{v'}$ and $q_{w'}$. We will finally show that the pair $(q_{v'},q_{w'})$ we found has indeed the given order-data.
 
    Each two consecutive kneading-sequences of \textbf{K} are distinct. Also, each ${\cal{K}}(c_{i})$ hits each plateau of $st{w} \circ st_{v}$ at most once, above its corresponding critical point. So, by lemma 2.4,  all kneading sequences of \textbf{K} are tight.  

    By Thurston's Theorem, these imply existence and uniqueness of a polynomial $P$ with kneading-data \textbf{K}, of shape (+,-,+,-) and conditions at the boundary $P(0)=0$ and $P(1)=0$ . A boundary anchored polynomial $P$ of degree 4, shape (+,-,+,-) and real distinct critical points $0<C_{1}<C_{2}<C_{3}<1$ is a composition of logistic maps if and only if $P(C_{1})=P(C_{3})$. Indeed, we know that the kneading sequences ${\cal{K}}(C_{1})={\cal{K}}(c_{1})$ and ${\cal{K}}(C_{3})={\cal{K}}(c_{3})$ are identical. Suppose $P(C_{1})<P(C_{3})$. Then the whole interval $[P(C_{1}),P(C_{3})]$ will have the same (bicritical) itinerary, as ${\cal{K}}(C_{1})={\cal{K}}(C_{3})$, so, after a finite number of iterations under $P$, it will all map to $C_{2}$, contradiction. So $P(C_{1})=P(C_{3})$, hence there exists a pair of quadratic maps such that $P=q_{w'} \circ q_{v'}$.

   The kneading data \textbf{K} determines the itinerary of the bicritical orbit and its order-data. So the polynomial map we found can only have the given order-data $(\sigma,\tau)$. } 
   
\vspace{3mm}

   Very similarly we can prove the equivalent statement for secondary intersections:

   \begin{thm} Let $(\sigma,\tau) \in S_{m+n}^{2}$ admissible joint order-data. There is a unique secondary intersection in $P^{Q}$ with this data and conversely.
   \end{thm}

\subsection{The correspondence of the \\  boundary points} \label{section:subsection9}
\indent

   Fix $(\sigma_{1},\tau_{1}) \in S_{n}^{2}$. The left bone $B^{ST}=B_{L}^{ST}(\sigma_{1},\tau_{1})$ in $P^{ST}$ with order-data $(\sigma_{1},\tau_{1})$ is as an algebraic curve in $P^{ST}=I_{2} \times I_{1} = I^{2}$. Its boundary consists of two points:
  
  $$ \delta B^{ST}=B^{ST} \cap \delta P^{ST}=B^{ST} \cap (I_{2} \times 
\{ 1 \})=\{ (v_{1},1),(v_{2},1) \} $$
with $v_{1} < v_{2}$.

   For any $(v,w)$, I will call $\Im_{ST}(x)(v,w)$ the itinerary of $x$ under $(st_{v},st_{w})$ and ${\bf K_{Q}}(v,w)$ the kneading-data of $st_{W} \circ st_{v}$. 

   The itineraries of the critical points $\gamma_{1}$ and $\gamma_{2}$ under $(st_{v_{1}},st_{1}$ and $(st_{v_{2}},st_{1})$ are respectively:\\

  $ \Im_{ST}(\gamma_{1})(v_{1},1) \neq \Im_{ST}(\gamma_{1})(v_{2},1) $\\
  
  $ \Im_{ST}(\gamma_{2})(v_{1},1) = \Im_{ST}(\gamma_{2})(v_{2},1) = 
                (\Gamma_{2},R_{1},L_{2},L_{1},L_{2},L_{1},...) = (\Gamma_{2},R_{1},\overline{L_{2},L_{1}})$ \\

   At any $(v,w) \in B^{ST}$, $\gamma_{1}$ has a periodic orbit $o_{1}$ of period $2n$ and order-data $(\sigma_{1},\tau_{1})$. At the two boundary points $(v_{1},1),(v_{2},1) \in \delta B^{ST}$, the orbit $o_{2}$ of $\gamma_{2}$ is also finite, although not periodic. 

    Statements in previous sections referred to primary or secondary intersections of bones. I will need some extensions of these statements to apply to boundary points of left bones in either parameter space. As we have noted, these boundary points are not bones crossings.\\

   We expect the boundary of the corresponding quadratic left bone $B^{Q}=B^{Q}(\sigma_{1},\tau_{1})$ to look similarly.

    \begin{thm} The boundary of $B^{Q}(\sigma_{1},\tau_{1})=B^{Q}$ consists of exactly two distinct points in $[0,1] \times \{1\} \subset \delta P^{Q}$. 
    \end{thm}

   \proof{ Consider the corresponding $ST$-left bone $B^{ST}(\sigma_{1},\tau_{1})$ and its boundary points $(v_{1},1)$ and $(v_{2},1)$. The maps $st_{v_{i}} \circ st_{1}$ have kneading-data ${\bf K_{ST}}(1,v_{i})$. For each $i$, the adjacent kneading-sequences are distinct.

   For each $i \in \{ 1,2 \}$, the pair of critical itineraries at $(1,v_{i})$ determines the respective kneading-data ${\bf K_{ST}}(1,v_{i})$. Note that $\Im_{ST}(\gamma_{1})(v_{1},1) \neq \Im_{ST}(\gamma_{1})(v_{2},1)$, so ${\bf K_{ST}}(1,v_{1}) \neq {\bf K_{ST}}(1,v_{2})$. The kneading-data also satisfies for each $i$ the conditions in the extended version of Thurston's theorem: the kneading sequences are finite and tight and  ${\cal}{K_{ST}}(1,v_{i})(c_{1})={\cal}{K_{ST}}(1,v_{i})(c_{3}) \neq {\cal}{K_{ST}}(1,v_{i})(c_{2})$. Hence for each $i$ there exists a point $(w'_{i},v'_{i}) \in P^{Q}$ such that $q_{v'_{i}} \circ q_{w'_{i}}$ has kneading-data ${\bf K_{Q}}(w'_{i},v'_{i})={\bf K_{ST}}(1,v_{i})$, and subsequently the same critical itineraries as $(st_{v_{i}},st_{1})$. In consequence: \\

   $ \Im_{Q}(\gamma_{1})(v'_{i},w'_{i})=\Im_{ST}(\gamma_{1})(v_{i},w_{i}) $ \\

   $ \Im_{Q}(\gamma_{2})(v'_{i},w'_{i}) = \Im_{ST}(\gamma_{2})(v_{i},1) = (\Gamma_{2},R_{1},\overline{L_{2},L_{1}})$ \\

   So clearly $(v'_{i},w'_{i})$ must be in the left bone $B_{L}^{Q}=B^{Q}$ in $P^{Q}$ corresponding to $B_{L}^{ST}=B^{ST}$ in $P^{ST}$. We also get that the itinerary of $\gamma_{2}$ under $(q_{v'_{i}},q_{w'_{i}})$ is $(\Gamma_{2},R_{1},\overline{L_{2},L_{1}})$. If $(v,w) \in [0,1]^{2}$ such that $vw>\frac{1}{16}$, then zero is a repeller for the composition $q_{w} \circ q_{v}$. This will be the case if we are situated on a left quadratic bone. So the only way for the itinerary of a point to stay indefinitely on $L_{1}$ and $L_{2}$ is for the point to map to zero after a number of iterates. To be consistent with the required itinerary, we need to have $(q_{v'_{i}} \circ q_{w'_{i}})(\gamma_{2}) = 0 $ and $q_{w'_{i}}(\gamma_{2})$ is $R$, so $q_{w'_{i}}(\gamma_{2})=1$, hence $w'_{i}=1$, for both $i=1$a d $i=2$.

   In conclusion: for the two points $(v_{1},1), (v_{2},1) \in \delta B^{ST}$ we found two points $(v'_{1},1), (v'_{2},1) \in \delta B^{Q}$ with the same corresponding kneading-data. The two points $(v_{1}',1)$ and $(v_{2}',1)$ we found in $\delta B^{Q}$ are the only two boundary points of $B^{Q}$. This follows almost immediately from Thurston's uniqueness.}

\subsection{A more complete description of  \\   bones in $ P^{ST}$ and $P^{Q}$} \label{section:subsection10}

\indent

  We plan to prove next: following the crossings along $B^{Q}=B^{Q}(\sigma_{1},\tau_{1}) \subset P^{Q}$, the combinatorics is same as at the crossings along the corresponding bone $B^{ST}=B^{ST}(\sigma_{1},\tau_{1}) \subset P^{ST}$. 

  We show first a combinatorial result concerning the order of occurrence of the primary and secondary intersections along a bone in $P^{ST}$ with fixed order-data $(\sigma_{1},\tau_{1})$. To fix our ideas, all proofs and results are developed for left bones $B^{ST}=B_{L}^{ST}$, hence we will omit writing the index $L$ unless it causes ambiguity.

   Fix a stunted left bone $B^{ST}=B^{ST}(\sigma_{1},\tau_{1})$ and slide $(v,w)$ along $B^{ST}$. Clearly, $\Im_{ST}(\gamma_{1})$ only changes at the primary intersection $(v_{0},w_{0})$. Therefore, $B_{*}^{ST}=B^{ST} \backslash \{(v_{0},w_{0})\}$ can be divided into two halves, each corresponding to a different itinerary of $\gamma_{1}$ under $(st_{v},st_{w})$; call $B_{-}^{ST}$ the left half, containing the boundary point $(v_{1},1) \in \delta B^{ST}$ and $B_{+}^{ST}$ the one containing $(v_{2},1) \in \delta B^{ST}$ (where $v_{1}<v_{2}$):

    $$ B^{ST} = B_{*}^{ST} \cup \{ (v_{0},w_{0}) \} = 
                B_{-}^{ST} \cup \{ (v_{0},w_{0}) \} \cup B_{+}^{ST} $$

   To fix our ideas, we look at $B_{-}^{ST}$; the results and their proofs should work symmetrically for $B_{+}^{ST}$. $B_{-}^{ST}$ is composed of a vertical segment and a horizontal one:

   $$ B_{-}^{ST}=\{ v_{1} \} \times [w_{0},1] \cup [v_{1},v_{0}] \times 
                 \{ w_{0} \} = B_{-,v}^{ST} \cup B_{-,h}^{ST} $$

\begin{figure}[h]
\begin{center}
\input{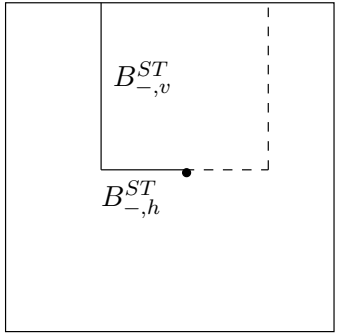}
\caption{\emph{ We divide the left half $B_{-}^{ST}$ of a left bone in $P^{ST}$ into a vertical segment $B_{-,v}^{ST}$ and a horizontal segment $B_{-,h}^{ST}$. All secondary intersections occur along $B_{-,v}^{ST}$. All points along the horizontal part are capture points.}}
\end{center}
\end{figure}

   We can now state our claim for this section in more precise terms:

   \begin{thm}  The secondary intersections occur along $B_{-,v}^{ST}$ in the strictly decreasing order of their itinerary $\Im_{ST}(\gamma_{2})$, as $w$ decreases from $1$ to $w_{0}$.
   \end{thm}

   \proof{ For a fixed $m \geq 1$, call ${\cal{D}}^{m}_{ST}$ the set of all parameters $(v,w)$ (secondary intersections and capture points) on $B_{-,v}^{ST}$ for which $\gamma_{2}$ maps to either $\gamma_{1}$ in $2m-1$ iterates or to $\gamma_{2}$ in $2m$ iterates. Call ${\cal{D}}_{ST}=\bigcup_{m \geq 1}{{\cal{D}}^{m}_{ST}}$ the \emph{distinguished} points on $B_{-,v}^{ST}$. Also call $\Im_{ST}^{m}(\gamma_{2})$ the itinerary $\Im_{ST}(\gamma_{2})$ truncated to the first $2m$ positions. 

   As $w$ decreases from $1$ to $w_{0}$, $\Im_{ST}^{m}(\gamma_{2})$ decreases (in the order inherited from the total order on infinite itineraries), with actual changes at all points in $\bigcup_{k \leq m}{{\cal{D}}^{k}_{ST}}$. Hence $\Im_{ST}(\gamma_{2})$ decreases, with changes at all points in ${\cal{D}}_{ST}$.

   Subsequently, $\Im_{ST}(\gamma_{2})$ decreases strictly on the set of distinguished points, in particular on the set of secondary intersections (see ~\cite{R} for details).}

\vspace{5mm}  

  \noindent \textbf{Remark.} The theorem makes it possible to identify the order of occurrence of the distinguished points (in particular of the secondary intersections) along $B_{-,v}^{ST}$ by looking at the itinerary of $\gamma_{2}$. From the construction of the stunted bones it is also easy to see that there are no secondary intersections on the horizontal segment of $B_{-,h}^{ST}$. In fact, all points of $B_{-,h}^{ST}$ are capture points and $\Im_{ST}(\gamma_{2})(v,w_{0})$ is constant for $v \in [v_{1},v_{0}]$.

    We move our focus now to the parameter space $P^{Q}$. The corresponding left bone $B^{Q}$ is a connected arc joining two boundary points $(v'_{1},1)$ and $(v'_{2},1)$ (with $v_{1}'<v_{2}'$) and having a unique primary intersection $(v'_{0},w'_{0})$. As before, the itinerary of $\gamma_{1}$ under $(q_{v'},q_{w'})$ changes only at $(v'_{0},w'_{0})$ as we move $(v',w')$ along $B^{Q}$. Hence we can divide $B^{Q}$ into two halves: left of $(v'_{0},w'_{0})$, containing $(v'_{1},1)$ and right of $(v'_{0},w'_{0})$, containing $(v'_{2},4)$.

  $$ B^{Q} = B_{-}^{Q} \cup \{ (v'_{0},w'_{0}) \} \cup B_{+}^{Q} $$
I will study the left half, comparatively with the vertical left half $B_{-,v}^{ST}$.

    We know that there is a bijective correspondence between secondary intersections along  $B_{-,v}^{ST}$ and $B_{-}^{Q}$ that associates to each intersection in $B_{-,v}^{ST}$ one with $\Im_{Q}(\gamma_{2})=\Im_{ST}(\gamma_{2})$ in $B_{-}^{Q}$ .We would like to prove that these secondary intersections occur on both $B_{-,v}^{ST}$ and $B_{-}^{Q}$ in the same decreasing order of $\Im(\gamma_{2})$, going from the boundary towards the primary intersection. In other words, we prove that the bijection is order preserving.

  Fix $m \geq 1$. Call $(l_{1},m_{1})$ the first distinguished point in $\bigcup_{k \leq m}{{\cal{D}}^{k}_{Q}}$ on $B_{-}^{Q}$ (from $(v'_{1},1)$ along the connected curve, with the regular order inherited by the order on $(0,1) \subset \Real$). 
   
  From theorem 2.15 we know that there is a corresponding distinguished point $(\alpha,\beta) \in \bigcup_{k \leq m}{{\cal{D}}^{k}_{ST}} \subset B_{-,v}^{ST}$ with the same critical itineraries :

  \begin{enumerate} 
  
    \item $ \Im_{ST}(\gamma_{1})(\alpha,\beta)=\Im_{Q}(\gamma_{1})(l_{1},m_{1})$ and 

    \item $ \Im_{ST}(\gamma_{2})(\alpha,\beta)=\Im_{Q}(\gamma_{2})(l_{1},m_{1})$ 
 
   \end{enumerate}
 
   \vspace{1mm}
 
   \noindent \textbf{Claim.} $(\alpha,\beta)$ is the first point to occur in $\bigcup_{k \leq m}{{\cal{D}}^{k}_{ST}}$ along $B_{-,v}^{ST}$.

   Suppose not. Then there exists a point $(v^{*},w^{*}) \in \bigcup_{k \leq m}{{\cal{D}}^{k}_{ST}}$ between the boundary point $(v_{1},1)$ and $(\alpha,\beta)$. We then have:\\

   $\Im_{ST}^{m}(\gamma_{2})(v_{1},1) > \Im_{ST}^{m}(\gamma_{2})(v^{*},w^{*}) > \Im_{ST}^{m}(\gamma_{2})(\alpha,\beta)$\\

   $\Im_{ST}^{m}(\gamma_{2})(v_{1},1)=\Im_{Q}^{m}(\gamma_{2})(v'_{1},1)$\\

   $\Im_{Q}^{m}(\gamma_{2})(l_{1},m_{1})=\Im_{ST}^{m}(\gamma_{2})(\alpha,\beta)$\\

   \noindent  The contradiction follows easily. (Note, for instance, that the conditions imply that the pair of critical itineraries at $(v_{1},1)$ has to be the same as the pair at a point right before $(\alpha,\beta)$).

   So the distinguished point in $(\alpha,\beta) \in B^{ST}_{-}$ with itinerary  $\Im_{ST}(\gamma_{2})(\alpha,\beta)=\Im_{Q}(\gamma_{2})(l_{1},m_{1})$ is the first to occur in $\bigcup_{k \leq m}{{\cal{D}}^{k}_{ST}}$. Continuing the procedure shows that the order of occurrence of all points in $\bigcup_{k \leq m}{{\cal{D}}^{k}_{ST}}$ along $B^{ST}_{-,v}$ is the same as the order of points in $\bigcup_{k \leq m}{{\cal{D}}^{k}_{Q}}$ along $B^{Q}_{-}$ (i.e. the decreasing order of the itinerary $\Im^{m}(\gamma_{2})$). We can state this as follows.

  \begin{thm}  For a fixed $m \geq 1$, going along $B^{ST}_{-,v}$ from $(v_{0},w_{0})$ to $(v_{1},1)$ and along $B^{Q}_{-}$ from $(v'_{0},w'_{0})$ to $(v'_{1},1)$, the itinerary $\Im^{m}(\gamma_{2})$ is monotonely increasing, with actual changes occurring at each distinguished point in $\bigcup_{k \leq m}{{\cal}{D}^{k}_{ST}}$ and  $\bigcup_{k \leq m}{{\cal}{D}^{k}_{Q}}$, respectectively. Hence the infinite itinerary $\Im(\gamma_{2})$ is monotonely increasing along  $B^{ST}_{-,v}$.
   \end{thm}

   \begin{thm} For a fixed $m \geq 1$, going along $[0,1] \times \{ 1 \} \subset\partial  P^{ST}$ and $[0,1] \times \{ 1 \} \subset  \partial P^{Q}$, the itinerary $\Im(\gamma_{2})=(\Gamma_{2},R_{1},\overline{L_{2},L_{1}})$ stays constant, but the itinerary $\Im^{m}(\gamma_{1})$ increases monotonically, with an actual change at each end-point of a bone of period $2k \leq 2m$.
   \end{thm}

\subsection{The big picture} \label{section:subsection11}
\vspace{5mm}

   \underline {\textbf{Overview of results}}: 

   \emph{Fix $m \geq 1$. Going along $B_{-}^{ST}$ from $(v_{0},w_{0})$ to $(v_{1},1)$ and along $B_{-}^{Q}$ from $(v'_{0},w'_{0})$ to $(v'_{1},1)$, the truncated itinerary $\Im^{m}(\gamma_{2})$ increases monotonically, with an actual increase at each crossing with a right bone. There is a one-to-one correspondence between the crossing points of bones of period at most $2m$ in the two families, correspondence that preserves the order of critical itineraries (i.e. of the joint order-data).}

   \emph{ Slide from left to right along the upper boundary of the two parameter spaces ($[0,1] \times \{ 1 \} \subset \partial P^{ST}$ and $[0,1] \times \{ 1 \} \subset \partial P^{Q}$). The itinerary $\Im(\gamma_{2})$ does not change, and the truncated itinerary $\Im^{m}(\gamma_{1})$ increases monotonely, with an actual change at each end-point of a left bone. There is a one-to-one correspondence between all boundary points of bones of period smaller than $2m$ in the two families, correspondence that preserves the order of the critical itineraries.}

\vspace{3mm}

   We want to restate the results in terms of kneading-data. In essence, we are looking to obtain in $P^{Q}$ a similar property to the following in $P^{ST}$ (see ~\cite{R}):

   We will use the following lemma:

   \begin{lemma} (a) Consider two arbitrary $(v'_{1},w'_{1}),(v'_{2},w'_{2}) \in B_{-}^{Q}$ and the itineraries $\Im_{Q}^{i}(\gamma_{2})$ of $\gamma_{2}$ under $q_{w'_{i}} \circ q_{v'_{i}}$, for $i=1,2$. If  $\Im_{Q}^{1}(\gamma_{2}) < \Im_{Q}^{2}(\gamma_{2})$ then the kneading data ${\bf K}(q_{w_{1}} \circ q_{v_{1}}) << {\bf K}(q_{w_{2}} \circ q_{v_{2}})$.

   (b) If $(v'_{1},1),(v'_{2},1) \in [0,1] \times \{ 1 \}$ are such that $\Im_{Q}^{1}(\gamma_{1}) < \Im_{Q}^{2}(\gamma_{1})$, then ${\bf K}(q_{1} \circ q_{v_{1}}) << {\bf K}(q_{1} \circ q_{v_{2}})$.
   \end{lemma}

\vspace{5mm}
   
   We restate two important conclusions in $P^{Q}$.

   \begin{thm} In the parameter space $P^{Q}$, the kneading-data of the maps $q_{w'} \circ q_{v'}$ increases along a left bone-arc from its primary intersection towards either boundary point and increases along the upper boundary interval $[0,1] \times \{1\} \in \partial P^{Q}$ from left to right (see picture). A symmetric statement holds for right bones and the right boundary interval.
   \end{thm}

   We know (see for example ~\cite{MT}) that the order of the kneading-data of two maps is preserved into the order of their topological entropies. Hence:

   \begin{thm} The topological entropy increases in $P^{Q}$ along each bone-arc from its primary intersection towards the boundary $\partial P^{Q}$ and along the boundary segments $[0,1] \times \{1\}$ and $\{1\} \times [0,1]$ towards the upper right corner (see picture).
   \end{thm}

\begin{figure}[h]
\begin{center}
\input{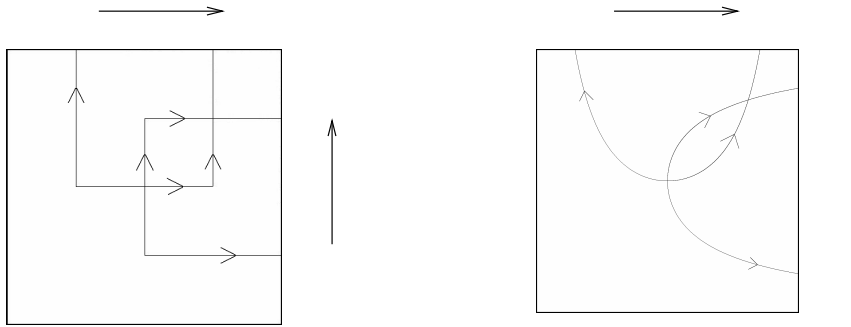}
\caption{\emph{ The arrows show the direction of increasing entropy along the bones and the boundary in $P^{ST}$ and $P^{Q}$. }}
\end{center}
\end{figure}

\vspace{5mm}

   We want to point out a few major consequences of our results, crucially important for later goals.

    We showed that every bone in $P^{Q}$ is composed of a bone-arc (that we called $B^{Q}$ in a previous section) and possible loop components. We will eventually rule out the existence of bone-loops. For the time being, a step towards this conclusion follows as a consequence of Thurston's uniqueness: \emph{for any arbitrary left bone in $P^{Q}$, the bone-arc $B^{Q}$ contains all possible post-critically finite kneading data (itineraries) admissible for the given bone.} In consequence, any loop component that the bone may have can not contain any post-critically finite points.

   \begin{defn} \emph{Fix $n \in \mathbb{N}$. We define the $n$-skeleton in either parameter space to be :}

    $S_{n}^{ST}$ = the union of all (left and right) bones $B_{2k}^{ST} \subset P^{ST}$ of  period $2k \leq 2n$, together with the boundary $\partial P^{ST}$;

    $S_{n}^{Q}$ = the union of all (left and right) bones $B_{2k}^{Q} \subset P^{Q}$ of  period $2k \leq 2n$, together with the boundary $\partial P^{Q}$.

    \emph{By a vertex of either skeleton we mean either an end-point of its bones or a (primary or secondary) intersection point.}
    \end{defn}

   \begin{thm} For any fixed $n \in \mathbb{N}$, there is a homeomorphism: 

   $$\eta_{n}:P^{ST} \; \longrightarrow \; P^{Q}$$ 

\noindent
which maps $S_{n}^{ST}$ onto $S_{n}^{Q}$, carrying $\partial P^{ST}$ to $\partial P^{Q}$, carrying bones to corresponding bones and and vertexes to vertexes with the same data.
   \end{thm}

   \proof{ We use the result that will be proved independently in the next two chapters: the bones in $P^{Q}$ are smooth ${\cal{C}}^{1}$ curves, intersecting transversally with each other and with the boundary. There are no bone loops in $P^{Q}$, so each bone is a smooth arc connecting two boundary points. Moreover, each such bone-arc contains all post-critically finite kneading-data existing on the corresponding bone in $P^{ST}$, in the same order of occurrence.

   The construction of the homeomorphism is topologically straightforward. Define $\eta_{n}$ on the set of vertexes by corresponding to each vertex in $S_{n}^{ST}$ the unique one in $S_{n}^{Q}$ with the same data. Along each bone, $\eta_{n}$ preserves the order of the vertexes. Hence we can extend it continuously to the intervals on the bones or boundary between each two vertexes, then to each skeleton-enclosed region. This can easily be done in such a way that the resulting continuous map $\eta_{n}:P^{ST} \; \longrightarrow \; P^{Q}$ is a homeomorphism. } 

\vspace{3mm}

\begin{figure}[h]
\begin{center}
\input{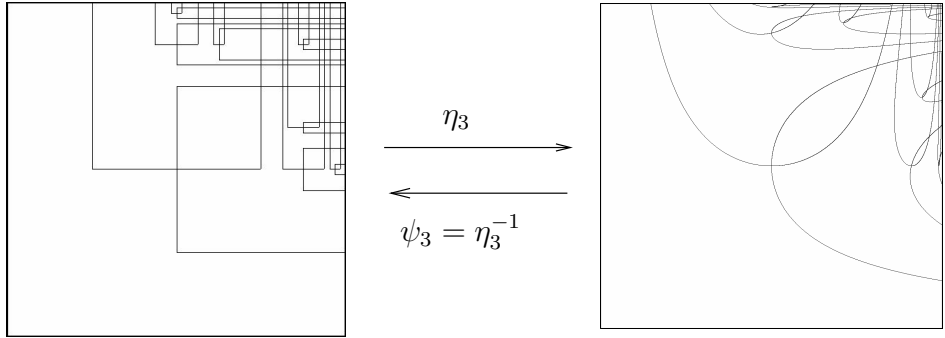}
\caption{\emph{ The $n$-skeletons define topological cell-complexes in both parameter spaces. The map $\eta_{n}$ is a homeomorphism between these complexes. The picture illustrates $n=3$.}}
\end{center}
\end{figure}

   We can associate to the $n$-skeleton in either parameter space a topological cell-structure as follows:

   $\bullet$ the 0-cells are points, more precisely the vertexes of the $n$-skeleton;

   $\bullet$ the 1-cells are the connected components of the bones obtained by deleting the vertexes, hence they are homeo to open intervals;

   $\bullet$ the 2-cells are the connected components of the complement of the $n$-skeleton in the respective parameter space, hence they are homeo to open discs.

   We will also use the closures of such cells, which are homeo to points, closed intervals and closed discs respectively.

\vspace{5mm}

   We call the resulting complexes: $P_{n}^{ST}$ in $P^{ST}$ and $P_{n}^{Q}$ in $P^{Q}$. The map $\eta_{n}:P_{n}^{ST} \; \longrightarrow \; P_{n}^{Q}$ is a homeomorphism of cell complexes, taking each cell in $P_{n}^{ST}$ to a corresponding cell in $P_{n}^{Q}$ by carrying vertexes to vertexes with the same entropy and edges to edges with the same interval of entropies.

\section{Hyperbolicity in $P^{Q}$}
\subsection{The mapping schema of a hyperbolic map} \label{section:subsection1}
\indent

    \begin{defn} Let $M$ be a finite disjoint union of copies of $\mathbb{C}$ and let $f:M \longrightarrow M$ be a proper holomorphic map of degree $\geq 2$ on each component of $M$. We say that $f$ is hyperbolic if every critical orbit converges to an attracting cycle.
   \end{defn}

   Let $f$ be a hyperbolic map as above. Let $W(f)$ be the union of the basins of attraction of all attracting cycles of $f$. $f$ carries each component $W_{\alpha} \subset W(f)$ onto a component $W_{\beta}$ by a map of degree $d_{\alpha} \geq 1$. Also let $W^{c}(f)$ be the union of all critical components $W_{\alpha} \subset W(f)$, that is of all $W_{\alpha}$ that contain critical points of $f$.

   We define the reduced mapping schema $\overline{S}(f)=(\mid S \mid ,F,w)$ associated to $f$ as the triplet made of:
 
   $\bullet$ a set of vertexes $\mid S \mid$, obtained by associating a vertex $\alpha$ to each critical component $W_{\alpha} \subset W^{c}(f)$;

   $\bullet$ a weight function $w: \mid S \mid \longrightarrow \mid S \mid$, defined as $w(\alpha)=$ the number of critical points of $f$ in $W_{\alpha}$;

   $\bullet$ a set of edges $F: \mid S \mid \longrightarrow \mid S \mid \; , \;F(\alpha)=\beta$, where $W_{\beta}$ is the image of $W_{\alpha}$ under the first return map to $W^{c}(f)$.

   The critical weight of $\overline{S}(f)$ is defined as $w(f)=\sum_{\alpha}{w(\alpha)}$

    All hyperbolic maps that interest us have reduced mapping schemata of critical weight 2, so we will only look at the cases that appear for $w=2$. For a more general analysis, see ~\cite{M1}.

   To a fixed mapping schema with $w=2$, we associate the universal polynomial model space ${\cal{P}}$. This will be the space of all maps $f$ from $\mathbb{C}_{1} \sqcup \mathbb{C}_{2}$ to itself such that the restriction of $f$ to each copy of $\mathbb{C}$ is a monic centered polynomial of degree 2. More precisely:\\

   $f(z)=z^{2}-a_{1}, \;$ for all $z \in \mathbb{C}_{1}$
   
   $f(z)=z^{2}-a_{2}, \;$ for all $z \in \mathbb{C}_{2}$\\
where $a_{1},a_{2} \in \mathbb{C}$.\\

We say that a map $f \in {\cal{P}}$ belongs to the connectedness locus ${\cal{C}}$ if its filled Julia set $K(f)$ intersects both $\mathbb{C}_{1}$ and $\mathbb{C}_{2}$ in a connected set. The hyperbolic connectedness locus ${\cal{H}} \subset {\cal{C}}$ is the open set of all $f \in {\cal{P}}$ for which the orbits of both critical points $0 \in \mathbb{C}_{1}$ and $0 \in \mathbb{C}_{2}$ converge to attracting periodic orbits.

   Such hyperbolic maps can be roughly classified into the three following types (see ~\cite{M3}):\\

  (1) \textbf{Bitransitive case}: $0 \in \mathbb{C}_{1}$ and $0 \in \mathbb{C}_{2}$ belong to $U_{1} \subset \mathbb{C}_{1}$ and $U_{2} \subset \mathbb{C}_{2}$ such that: $U_{1}$ is mapped to $U_{2}$ under $q_{1}$ iterates of $f$ and $U_{2}$ is mapped to $U_{1}$ under $q_{2}$ iterates.

\begin{figure}[h]
\begin{center}
\input{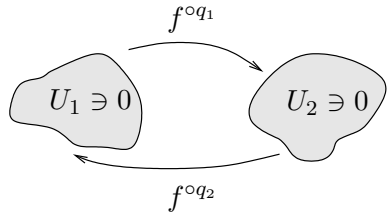}\\
\caption{\emph{ The behavior of a bitransitive hyperbolic map.}} 
\end{center}
\end{figure}

   (2) \textbf{Capture case}: $0 \in U_{1} \subset \mathbb{C}_{1}$ and $0 \in U_{2} \subset \mathbb{C}_{2}$ such that $U_{1}$ is periodic and $U_{2}$ is not, but some forward image of $U_{2}$ coincides with $U_{1}$. Also its symmetric case.

\begin{figure}[h]
\begin{center}
\input{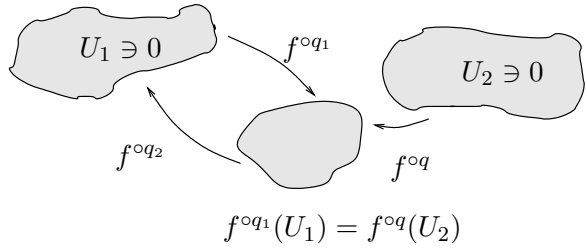}\\
\caption{\emph{The behavior of a map in the capture case.}}
\end{center}
\end{figure}

   (3) \textbf{Disjoint periodic sinks}: $0 \in U_{1}$ and $0 \in U_{2}$, where $U_{1}$ and $U_{2}$ are periodic of periods $q_{1}$ and $q_{2}$, but no forward image of $U_{1}$ coincides with $U_{2}$ and vice-versa.\\

\begin{figure}[h]
\begin{center}
\input{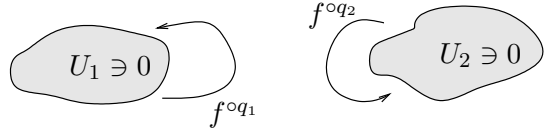}\\
\caption{\emph{ The behavior of a map in the disjoint sinks case.}}
\end{center}
\end{figure}

   For maps $f \in {\cal{H}}$, we may consider their reduced mapping schemata $\overline{S}(f)$. These schemata will all have critical weight 2, but not all are isomorphic (see figure 14). However, all maps in each connected component of ${\cal{H}}$ clearly have isomorphic schemata. Furthermore, by theorem 4.1 in ~\cite{M1}:

   \begin{thm} If $H_{\alpha} \subset {\cal{C}}$ is a hyperbolic component of ${\cal{H}}$ with maps having reduced schemata isomorphic to $S$, then $H_{\alpha}$ is diffeomorphic to a model space $B(S)$. In particular, any two hyperbolic components $H_{\alpha}$ and $H_{\beta}$ with schemata isomorphic to $S$ are diffeomorphic. Moreover, each $H_{\alpha}$ contains a unique post-critically finite map $f_{\alpha}$, called its center.
   \end{thm}

\begin{figure}[h]
\begin{center}
\input{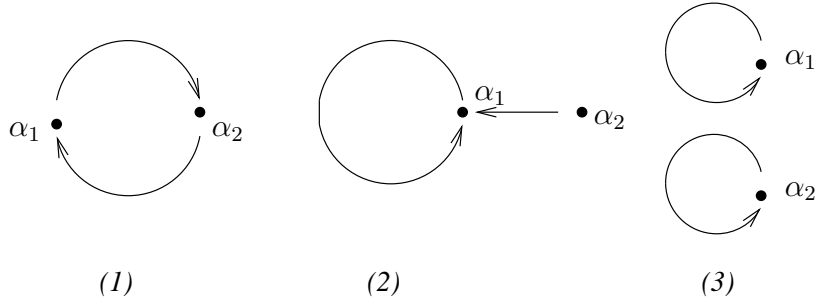}\\
\caption{\emph{ (1) Bitransitive case: $\mid S \mid =\{ \alpha_{1},\alpha_{2} \}, \quad F(\alpha_{1})=\alpha_{2}, \; F(\alpha_{2})=\alpha_{1}, \quad \omega(\alpha_{1})=\omega(\alpha_{2})=1$. (2) Capture case: $\mid S \mid =\{ \alpha_{1},\alpha_{2} \}, \quad F(\alpha_{1})=\alpha_{1}, \; F(\alpha_{2})=\alpha_{1}, \quad \omega(\alpha_{1})=\omega(\alpha_{2})=1$. (3) Disjoint sinks case: $\mid S \mid =\{ \alpha_{1},\alpha_{2} \}, \quad F(\alpha_{1})=\alpha_{1}, \; F(\alpha_{2})=\alpha_{2}, \quad \omega(\alpha_{1})=\omega(\alpha_{2})=1$}}
\end{center}
\end{figure}

   \begin{defn} A real form of the mapping schema $S$ is an antiholomorphic involution $\rho: \mathbb{C}_{1} \sqcup \mathbb{C}_{2} \longrightarrow \mathbb{C}_{1} \sqcup \mathbb{C}_{2}$ which commutes with the special map $f_{0}^{S}: \mathbb{C}_{1} \sqcup \mathbb{C}_{2} \longrightarrow  \mathbb{C}_{1} \sqcup \mathbb{C}_{2}, \;\, f_{0}^{S}(z)=z^{2}$. The collection of maps $f \in {\cal{P}}$ that commute with $\rho$ is an affine space ${\cal{P}}_{\mathbb{R}}(\rho)$, which we call the real form of ${\cal{P}}$ associated with $\rho$. We also define the corresponding real connectedness locus and the real hyperbolic locus as:\\

   ${\cal{C}}_{\mathbb{R}}(\rho)={\cal{C}} \cap {\cal{P}}_{\mathbb{R}}(\rho)$

   ${\cal{H}}_{\mathbb{R}}(\rho)={\cal{H}} \cap {\cal{P}}_{\mathbb{R}}(\rho)$
   \end{defn}

   For each mapping schema of weight 2, there are exactly two real forms. The form $\rho_{0}(z)=\overline{z}$ corresponds to the space ${\cal{P}}_{\mathbb{R}}(\rho_{0})$ of real polynomials in ${\cal{P}}$. If we restate theorem 6.4 of ~\cite{M1} in our particular case, we obtain:

   \begin{thm} Any hyperbolic component in ${\cal{C}}_{\mathbb{R}}={\cal{C}}_{\mathbb{R}}(\rho_{0}) \subset {\cal{P}}_{\mathbb{R}}(\rho_{0})$ is a topological 2-cell with a unique ``center point'' and is real analytically homeomorphic to a space of Blaschke products $\beta_{\mathbb{R}}(S,\rho_{0})$.
   \end{thm}

   In other words, all hyperbolic components with the same schemata in  ${\cal{C}}_{\mathbb{R}}$ are diffeomorphic to each other. For example, all bitransitive components are diffeo to the principal component centered at:

    $$f_{0}^{S}: \mathbb{C}_{1} \sqcup \mathbb{C}_{2} \longrightarrow  \mathbb{C}_{1} \sqcup \mathbb{C}_{2}, \;\, f_{0}^{S}(z)=z^{2}$$

    For a detailed characterization of the construction and properties of the suitable Blaschke-products model spaces, see ~\cite{M1}.
\subsection{ Hyperbolic components in $P^{Q}$} \label{section:subsection2}
\indent

    Let us return to our space, containing real quartic polynomials that are compositions $q_{w} \circ q_{v}$ of logistic maps.

   Let $\mathbb{C}_{1}$ and $\mathbb{C}_{2}$ be two copies of the complex plane and consider $q_{v}:\mathbb{C}_{1} \longrightarrow \mathbb{C}_{2}$ and $q_{w}:\mathbb{C}_{2} \longrightarrow \mathbb{C}_{1}$ the complex extensions of two fixed logistic maps of the interval. We define a new map: $q^{v}_{w}: \mathbb{C}_{1} \sqcup \mathbb{C}_{2} \longrightarrow \mathbb{C}_{1} \sqcup \mathbb{C}_{2}$, acting as $q_{v}$ on $\mathbb{C}_{1}$ and as $q_{w}$ on $\mathbb{C}_{2}$.

   Let $W(q^{v}_{w}) \subset \mathbb{C}_{1} \sqcup \mathbb{C}_{2}$ be the open set consisting of all complex numbers in $\mathbb{C}_{1}$ and $\mathbb{C}_{2}$ whose forward orbit under $q^{v}_{w}$ converges to an attracting periodic orbit of $q^{v}_{w}$.

  Under iteration of $q^{v}_{w}$, each component of $W(q^{v}_{w})$ is mapped onto a component of $W(q^{v}_{w})$. As before, we will say that $q^{v}_{w}$ is hyperbolic if both $\gamma_{1} \in I_{1} \subset \mathbb{C}_{1}$ and $\gamma_{2} \in I_{2} \subset \mathbb{C}_{2}$ are contained in $W(q^{v}_{w})$.

   It would be convenient to find a correspondence between our family of pairs of real quadratic maps, parametrized by $(v,w) \in P^{Q}$ and the family of degree $2$ normal polynomials. It can be shown that each map $q_{w} \circ q_{v}: \mathbb{C}_{1} \longrightarrow \mathbb{C}_{2}$ is conjugated by a complex affine map $L$ to a composition of maps $z \longrightarrow z^{2}-a_{1}$ and $z \longrightarrow z^{2}-a_{2}$. Moreover, the correspondence $(v,w) \rightarrow (a_{1},a_{2})$ is ``nice'' enough to permit us to carry over to $P^{Q}$ properties we have in the space of normal forms. More precisely:

   \begin{thm} Let $U$ be the subset of $P^{Q}$ consisting of pairs $(v,w)$ with $vw > \frac{1}{16}$. For each such pair $(v,w) \in U$ there is a unique pair $(A,B) \in \mathbb{R}^{2}$ such that $q_{w} \circ q_{v}$ is linearly conjugate to $z \longrightarrow z^{4}+Az^{2}+B$; there also exists a unique pair $(a_{1},a_{2}) \in \mathbb{R}^{2}$ so that $q_{w} \circ q_{v}$ is linearly conjugate to the composition of $z \longrightarrow z^{2}-a_{1}$ and $z \longrightarrow z^{2}-a_{2}$.

   Furthermore, recall that the connectedness locus ${\cal{C}}_{\mathbb{R}} \subset \mathbb{R}^{2}$ is the subset of parametes  $(a_{1},a_{2}) \in \mathbb{R}^{2}$  for which the complex critical points of  $(z^{2}-a_{1})^{2}-a_{2}$  have bounded orbits. The correspondence described above:

   $$ \Xi:\, U \longrightarrow {\cal{C}}_{\mathbb{R}}$$
   $$ \Xi(\lambda,\mu)=(a_{1},a_{2})$$
is a bijective diffeomorphism.
   \end{thm}

    \proof{ Each $q_{w} \circ q_{v}$ with $(v,w) \in P^{Q}$ is conjugated by an affine map $L(z)=-\frac{8}{\sqrt[3]{v^{2}w}}z+\frac{1}{2}$ to a composition of the two monic centered quadratic complex maps: $z \longrightarrow \zeta=z^{2}-a_{1}(\lambda,\mu)$ and $\zeta \longrightarrow z=w^{2}-a_{2}(\lambda,\mu)$. The correspondence:

     $$\Phi:U \longrightarrow \mathbb{R}^{2}, \; \Phi(v,w)=(A,B)$$
is a diffeomorphism onto its image, where the image $\Xi(U)$ is exactly the real connectedness locus ${\cal{C}}_{\mathbb{R}}$ in ${\cal{P}}_{\mathbb{R}}$.}

\vspace{3mm}

   \noindent \textbf{Remarks. } (1) The region $P^{Q} \,\backslash \,U = \{ (v,w) \, \slash \, vw <\frac{1}{16} \}$ is itself a hyperbolic component of $P^{Q}$, whose maps have all critical points attracted to zero. The map $\Xi$ folds this region and the principal component centered at $(v,w)=(\frac{}{2},\frac{1}{2}) \in P^{Q}$ onto the same component in ${\cal{C}}_{\mathbb{R}}$.

   (2) \emph{All bones in $P^{Q}$ are contained in $U$.} Indeed, suppose there is a $(v,w)$ on a bone such that $(v,w) \notin U$. The fixed origin is not repelling for the map $q_{w} \circ q_{v}$ with negative Schwarzian derivative, so it attracts all critical points, hence $(v,w)$ can't be on a bone, contradiction.

\vspace{3mm}

\begin{figure}[h]
\begin{center}
\input{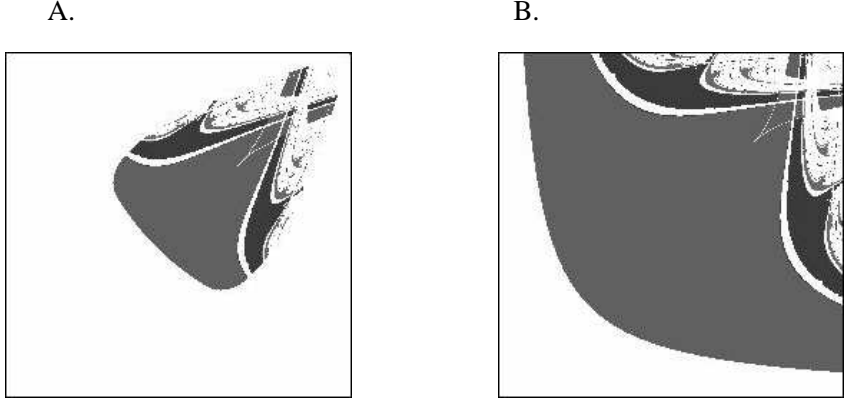}\\
\caption{\emph{ A. Hyperbolic components in $\mathbb{R}^{2}$ for the classical family of pairs quadratic monic centered maps. The picture shows the parameter window $(a_{1},a_{2}) \in [-2,2] \times [-2,2]$. B. Hyperbolic components in $U \in P^{Q}$. The principal component in both cases is visible as the large central shaded region.}}
\end{center}
\end{figure}

    We use the results in the previous sections to give the needed description of the hyperbolic components in our original parameter space $P^{Q}$. Hyperbolic components within each class (bitransitive, capture and disjoint sinks) are diffeomorphic to each other. The center points in each case will be respectively a primary intersection, a capture point or a secondary intersection.

   \begin{thm} Each hyperbolic component in $U \in P^{Q}$ is a topological 2-cell which contains a unique post-critically finite point, called its center. Moreover, every bone that intersects such a component does it along a simple arc passing through the center. Subsequently, there could be either one bone crossing the component through its center (capture case) or a pair of left-right bones intersecting transversally at the center point (bitransitive and disjoint sinks cases).
   \end{thm}

\subsection{ Density of hyperbolicity in $P^{Q}$} \label{section:subsection3}
\indent

   We aim to prove the following main result:

   \begin{thm} Hyperbolicity is dense in the parameter space $P^{Q}$.
   \end{thm}

\vspace{3mm}

\noindent \textbf{Remark.} The theorem is a modification of the more general Fatou conjecture (see ~\cite{KSvS}). The reference gives a proof that makes use of the following \emph{Rigidity Theorem}, that we will also be used to prove theorem 3.7.\\

\noindent \textbf{Rigidity Theorem.} \emph{Let $f$ and $f'$ be two polynomials with real coefficients, real non-degenerate critical points, connected Julia set and no neutral periodic points. If $f$ and $f'$ are topologically conjugate as dynamical systems on the real line $\mathbb{R}$, then they are quasiconformally conjugate as dynamical systems on the complex plane $\mathbb{C}$.}

\vspace{3mm}

  \proof{ We define the family ${\cal{S}}_{4}$ as the set of complex polynomials $Q: \, \mathbb{C} \, \rightarrow \, \mathbb{C}$ of degree 4, ``boundary anchored'' (i.e. $Q(0)=Q(1)=0$) and such that $Q(z)=Q(1-z)$, for all $z \in \mathbb{C}$.\\

   Consider $X_{s}$ to be the subset of maps in ${\cal{S}}_{4}$ with the following properties:

  $\bullet$ They have real coefficients.
  
  $\bullet$ Their three critical points are real and nondegenerate.

  $\bullet$ all critical points and values are in [0,1]. Hence their Julia sets are connected (see for example theorem 17.3 in ~\cite{M4}).

  $\bullet$ The boundary $\{0,1\}$ is repelling.\\

\begin{figure}[h]
\begin{center}
\input{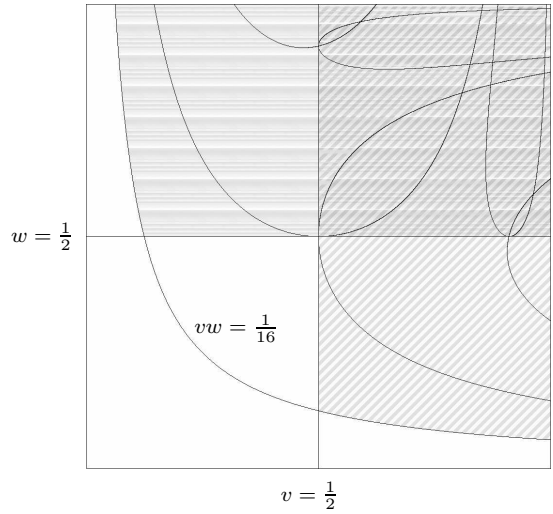}\\
\caption{\emph{ All maps in $\{ vw <\frac{1}{2} \}$ and in $\{ v <\frac{1}{2}, \, w <\frac{1}{2} \}$ are hyperbolic. Hyperbolic maps are dense in $\{ vw >\frac{1}{16}, \, v \geq \frac{1}{2} \}$ (slant shaded). By symmetry, they are dense in $\{ vw >\frac{1}{16}, \, w \geq \frac{1}{2} \}$ (horizontaly shaded). The region $vw >\frac{1}{16}$ contains all left and right bones.}}
\end{center}
\end{figure}

\noindent
In other words:

$$X_{s}=\{ q_{w} \circ q_{v}, \: \text{where} \; (v,w) \in P^{Q}, \, v \geq \frac{1}{2}, \, vw >\frac{1}{16} \}$$

\noindent
Indeed, recall that the three complex critical points of an arbitrary $P \in P^{Q}$ are $C_{1}$, $C_{2}=\frac{1}{2}$ and $C_{3}=-C_{1}$. An equivalent condition to $C_{1} \in \mathbb{R}$ is that:

   $$ q_{v}(\frac{1}{2}) \geq \frac{1}{2} \; \Leftrightarrow \; v \geq \frac{1}{2}$$

   We \textbf{claim} that hyperbolic polynomials are dense in $X_{s}$. Then the proof of 3.7 follows relatively easily. Indeed, the claim implies directly density of hyperbolicity in the region in $P^{Q}$ where $vw > \frac{1}{16}$ and $v \geq \frac{1}{2}$. By the symmetry property $(2)$, the result follows in the region where $vw > \frac{1}{16}$ and $w \geq \frac{1}{2}$. In the regions $\{ vw>\frac{1}{16}, \, v<\frac{1}{2},\, w<\frac{1}{2} \}$ and $\{ vw <\frac{1}{2} \}$ the proof is trivial: if $vw < \frac{1}{16}$ then all three critical orbits of $q_{w} \circ q_{v}$ converge to zero, while if $v<\frac{1}{2}$, $w<\frac{1}{2}$ and $vw >\frac{1}{16}$ then all critical orbits converge to a point in $(0,\frac{1}{2})$.}

\vspace{3mm}

   Next, we aim to prove density of hyperbolicity in $X_{s}$. \\

    \begin{lemma} Consider $P \in X_{s}$ with one parabolic cycle $\{ z_{1},..., z_{m} \}$. We can approximate $P$ by a polynomial $S \in X_{s}$ for which the cycle is attracting.
   \end{lemma}

\noindent \textbf{Sketch of proof:}  Fix $P \in X_{s}$ as above.

   It is fairly easy to show the existence of a polynomial $Q:\mathbb{C} \, \longrightarrow \, \mathbb{C}$ with real coefficients and the following properties (see ~\cite{R}):

   $\bullet \quad Q(z)=Q(1-z), \; \forall \, z \in \mathbb{C}$

   $\bullet \quad Q(z_{j})=0, \; \forall \, j=\overline{1,m}$

   $\bullet \quad Q(0)=Q(1)=0$

   $\bullet \quad Q'(x)=0$ when $P'(x)=0$

   $\bullet \quad \sum {\displaystyle \frac{Q'(z_{j})}{P'(z_{j})}}<0$

   Consider the new polynomial $R=P+\epsilon Q$. For small real values of $\epsilon$, $R$ perturbes the neutral cycle of $P$ to an attracting cycle:

    $$ \sum{ \log \mid R'(z_{j}) \mid} = \sum{ \log \mid P'(z_{j}) \mid} + \sum{ \log \mid 1+ \epsilon \frac{Q'(z_{j})}{P'(z_{j})} \mid } = $$
   $$ =\epsilon \sum{\frac{Q'(z_{j})}{P'(z_{j})}} + {\cal}{o}(\epsilon^{2}) < 0 $$

\vspace{3mm}

   For small enough values of $\epsilon$, $R$ has the following properties:

   $\bullet$ the parabolic cycle of $P$ is attracting for $R$;

   $\bullet$ the attracting/repelling cycles of $P$ change to attracting/repelling cycles for $R$ (hence $\{ 0 \}$ remains a repelling fixed boundary point for $R$);

   $\bullet$ $R(z)=R(1-z),\; \forall \, z \in \mathbb{C}$ and $R(0)=R(1)=0$, hence $R \in {\cal{S}}^{4}$;

   $\bullet$ R has real coefficients;

   $\bullet$ the critical points of $R$ are the same as the critical points of $P$, hence they are real, nondegenerate; all critical points and values are contained in $[0,1]$, hence the Julia set $J(R)$ is connected;\\

   However, in order to satisfy all required conditions, $Q$ (hence $R$) may have degree larger than 4. We use the Straightening Theorem to obtain a degree 4 polynomial $S \in X_{s}$ with the same behavior as $R$ (see for example ~\cite{CG}or ~\cite{R}). \hfill $\Box$

\vspace{5mm}

   For every $Q \in {\cal{S}}_{4}$, let $\tau(Q)$ be the number of critical points contained in the attracting basin of a hyperbolic attracting cycle of $Q$. Define: 

  $$ X'_{s} = \{ Q \in X_{s} \; / \; \tau(Q) \: \text{has a local maximum at} \: Q \}$$

   As $\tau$ is uniformly bounded above, $X'_{s}$ is dense in $X_{s}$. Moreover, $\tau$ is locally constant at any $P \in X'_{s}$, hence we have the following:

   \begin{prop} $X_{s}'$ is open and dense in $X_{s}$.
   \end{prop}

   \begin{prop} No map in $X'_{s}$ has a neutral cycle.
   \end{prop}

   \proof {Consider $P \in X'_{s}$ and $Q$ given by the lemma. By making the perturbation small enough, we can arrange that the other hyperbolic attractors of $P$ do not disappear. Moreover, we can also make sure that the critical points that were attracted to the attracting cycles remain so under the perturbation. 

   On the other hand, each attracting cycle attracts at least one critical point. Hence introducing a new attractor by perturbing $P$ to $Q$ will change $\tau$ as :

   $$\tau(Q) \geq \tau(P) + 1$$
contradiction with the local maximality of $\tau$ at $P$. }

   We finish by giving a reduced statement, from which theorem 3.7 follows now immediately. The proof is detailed in section 3.4.

   \begin{thm} Hyperbolic polynomials are dense in $X_{s}'$.
   \end{thm}

\subsection{ A reduced density result} \label{section:subsection4}
\indent

   Recall that two points $z_{1}$ and $z_{2}$ are in the same foliated equivalence class of a map $f$ if their grand orbits under $f$ have the same closure. For a fixed $f$, we denote by $n_{ac}$ the number of foliated equivalence classes of acyclic critical points in the Fatou set of $f$.
   By ~\cite{MS}, the complex dimension of the Teichmuller space of a map $f:\mathbb{C} \, \rightarrow \, \mathbb{C}$ is given by:\\

    $ \text{dim}(Teich(f)) = n_{ac}+n_{hr}+n_{lf}+n_{p} $,  where:\\

   $ n_{ac}$ = $\#$ of foliated equivalence classes of acyclic critical points in the Fatou set $F(f)$;\\

   $ n_{hr}$ = $\#$ of Herman rings of $f$;\\

   $ n_{lf}$ = $\#$ invariant line fields;\\

   $ n_{p}$  = $\#$ parabolic cycles.\\

   If $P \in X_{s}'$, $P$ has no Herman rings and no Siegel discs. By ~\cite{KSvS} and ~\cite{S}, $P$ does not support an invariant line field in its Julia set. We also proved in lemma  3.8 that $P$ does not have any parabolic basins. So all connected components of its Fatou set are attracting basins. Hence:

    $$ n_{hr} = n_{lf} = n_{p} = 0 \; \Rightarrow \; \text{dim}(Teich(P)) = n_{ac} $$

   Hence the set:\\

   $ QC(P) = \{ Q \in {\cal}{S}_{4} \; / \; Q \; \text{quasiconformally conjugate to} P \}$ \\

\noindent
is covered by countably many complex submanifolds of dimension $n_{ac}$. Subsequently, the set:\\

   $QC^{\mathbb{R}}(P)=QC(P) \cap X_{s}$\\

\noindent
is covered by countably many embedded real analytic submanifolds of $X_{s}$ with real dimension $n_{ac}$.

   We will also use the following(~\cite{dMvS}, pp 93):

   \begin{defn} \emph{If the 3-modal maps $P,Q: [0,1] \, \rightarrow \, [0,1]$ are such that}

    $h_{P}^{Q}: \; {\bigcup}_{n,i}{P^{n}(c_{i}(P))} \; \rightarrow \; {\bigcup}_{n,i}{Q^{n}(c_{i}(Q))} \; i=\overline{1,2,3}$

\noindent
\emph{defined by} :

   $h_{P}^{Q}(P^{n}(c_{i}(P)))=Q^{n}(c_{i}(Q)), \; \forall i=\overline{1,2,3} , \; \forall n \in \mathbb{N}$

\noindent 
\emph{is an order-preserving bijection, then we say that $P$ and $Q$ are combinatorially equivalent as 3-modal maps of the interval.}
   \end{defn}

   The relationship between combinatorial equivalence and topological conjugacy in our space $X_{s}$ can be described by the following theorem  (~\cite{dMvS}):\\

   \begin{thm}  Call ${\cal{F}}$ the family of maps $f$ of the interval satisfying the following:

   (1) they are of class ${\cal{C}}^{3}$;

   (2) they have nonflat critical points ( i.e. $D^{2}f(c) \neq 0, \; \forall c$ such that $Df(c)=0$ );
 
   (3) they have negative Schwartzian derivative: $Sf < 0$;

   (4) the boundary of the interval is repelling (in other words $\mid Df(x) \mid >1, \;$ if $x \in \{ 0,1 \}$);

   (5) they have no one-sided periodic attractors.

   Two maps $f,g \in {\cal{F}}$ are topologically conjugate ($f \stackrel{top}{\sim}_{_{\mathbb{R}}} g$) if and only if they are combinatorially equivalent ($f \stackrel{c.e.}{\sim}_{_{\mathbb{R}}} g$).
   \end{thm}

   \noindent \textbf{Remark.} If $P$ and $Q$ are maps in $X_{s}'$ restricted to the interval $[0,1]$, then both the conditions of theorem 3.13 and the Rigidity Theorem are satisfied, hence we have the following implications:

   $$P \stackrel{c.e.}{\sim}_{_{\mathbb{R}}} Q \; \Leftrightarrow \; P \stackrel{top}{\sim}_{_{\mathbb{R}}} Q \; \Rightarrow P \stackrel{qc.}{\sim}_{_{\mathbb{C}}} Q$$

   \vspace{5mm}

  \noindent \textbf{Proof of theorem ...}  Fix $P \in X_{s}'$.

   We think of $\; {\cal{S}}_{4} \subset \mathbb{C}^{2} \;$ and we consider the three holomorphic functions $c_{i} : {\cal{U}} \; \rightarrow \; \mathbb{C}, \; i=\overline{1,2,3}$ that give the three critical points of each map  $Q \in {\cal{U}}$. By taking ${\cal{B}} \subset {\cal{U}} \subset {\cal{S}}^{4}$ to be a small ball around $P$, we can arrange  to have $c_{1}(Q) < c_{2}(Q) < c_{3}(Q)=-c_{1}(Q)$, for any $Q \in {\cal{B}} \cap X_{s}$. Take ${\cal{B}}$ small enough for $\tau$ to be constant: $\tau=\tau(Q), \; \forall Q \in {\cal{B}} \cap X_{s}$ (recall $\tau$ is locally constant at each $P \in X_{s}'$).\\

   We want to prove (by contradiction) that ${\cal{B}} \cap X_{s}$ contains hyperbolic maps. Suppose the maps in ${\cal{U}} \cap X_{s}$ are not hyperbolic, hence $\tau<3$. There are two cases that remain for analysis:\\

   (1) $\tau=1$ (only $C_{2}$ is attracted) or $\tau=2$ (only $C_{1}$ and $C_{3}$ are attracted). Either way, there is only one foliated equivalent class of critical points in the Fatou set, hence $n_{ac} \leq 1$ (note that the critical points are not necessarily acyclic). Hence $QC^{\mathbb{R}}(Q)$ is in this case at most a countable union of lines in $X_{s}$, for any $Q \in {\cal{B}} \cap X_{s}$.\\

   (2) $\tau=0$ (no critical points are attracted). Hence $n_{ac}=0$, so $QC^{\mathbb{R}}(Q)$ is a countable union of points in $X_{s}$, for any $Q \in {\cal{B}} \cap X_{s}$.\\

   \textbf{A.} Suppose first there are no bones crossing the neighbourhood ${\cal{B}}$.\\

   If there are no other ``critical relations'' in ${\cal{B}}$ (i.e. there are no $m,n \in \mathbb{N}$ such that $Q^{m}(c_{1}(Q))=Q^{n}(c_{2}(Q))$ for some $Q \in {\cal{B}}$), then for any arbitrary $Q \in {\cal{B}}$ the map $h_{P}^{Q}$ defined in 3.12 is order preserving.(Note that we do not consider $Q(c_{1}(Q))=Q(c_{3}(Q))$ a critical relation.) Indeed: Suppose that $h$ reverses the order of two elements:

   $\quad P^{k}(c_{i}(P))<P^{l}(c_{j}(P))$ and

   $\quad Q^{k}(c_{i}(Q))>Q^{l}(c_{j}(Q))$  

   By continuity, there exists a $T \in {\cal{B}}$ such that:

   $\quad T^{k}(c_{i}(T))=T^{l}(c_{j}(T))$, contradiction.\\

   Since $h_{P}^{Q}$ is order-preserving for any $Q \in {\cal{B}} \cap X_{s}$, it follows that $P$ is combinatorially equivalent to any $Q \in {\cal{B}} \cap X_{s}$, hence $P$ is quasiconformally conjugate to any $Q \in {\cal{B}} \cap X_{s}$. This contradicts the fact that $QC^{\mathbb{R}}(P)$ is at most a union of countably many lines in $X_{s}$.

\vspace{3mm}

   Clearly, the ``no critical relations'' condition applies in the case $\tau=1$ or $\tau=2$. 

   If $\tau=0$, it could happen that all neibourhoods of $\tau$, arbitrarily small, contain critical relations. In other words, there exists a map $R$ arbitrarily close to $P$ that has a critical relation, say $R^{m}(c_{1}(R))=R^{n}(c_{2}(R))$.  

   Consider $\Sigma=\{ Q \in {\cal{B}} \cap X_{s} \, \slash \, Q^{m}(c_{1}(Q))=Q^{n}(c_{2}(Q)) \}$. This is a 1-dim curve in ${\cal{B}} \cap X_{s}$. There clearly are no other critical relations on $\Sigma$, hence the map $h_{R}^{Q}$ is order-preserving for any $Q \in \Sigma$. Subsequently, all maps in $\Sigma$ are combinatorially equivalent to $R$, hence quasiconformally conjugate to $R$. This contradicts the fact that $QC^{\mathbb{R}}(R)$ is a collection of countably many points in $X_{s}$, as $\tau=0$. \\

  \textbf{B.} If ${\cal{B}} \cap X_{s}$ is crossed by a bone $B$, let $R \in B \cap {\cal{B}} \cap X_{s}$.

   Bones can't accumulate at $R$, or $R$ would be hyperbolic. So there exists a neighbourhood ${\cal{V}}$ of $R$, ${\cal{V}} \subset {\cal{B}} \cap X_{s}$ that intersects no other bones than $B$. Take $S \in {\cal{V}} \backslash B$ and take ${\cal{W}}$ a neighbourhood of $S$ in ${\cal{V}} \backslash B$. Then the argument at \textbf{A.} applies for ${\cal{W}}$ and leads us to a contradiction.  \hfill $\Box$ 

\vspace{3mm}

   The proof of theorem 3.7 is now finished.   
\section{Topological properties of the Q-bones}
\subsection{Smoothness of the Q-bones} \label{section:subsection1}

   As we have stated before, a bone in $P^{Q}$ is an algebraic variety with two boundary points in $\partial P^{Q}$. As far as we presently know, the bone curves may not even be connected. We will rule this out in chapter 4.3, where we show independently that a bone can't contain any loops. For now, we dedicate this paragraph to proving that:

   \begin{thm} The bones are smooth ${\cal{C}}^{1}$ curves that intersect transversally.
   \end{thm}

   Recall that we use the notations $B_{L,2n}^{Q}$ and $B_{R,2n}^{Q}$ for a left/right bone in $P^{Q}$ of period $2n$ and given order-data. Fix an arbitrary point $p_{0}=(v_{0},w_{0})$ on a left bone $B^{Q}_{L,2n}$. We want to show that $B^{Q}_{L,2n}$ is smooth at $p_{0}=(v_{0},w_{0})$.

   For the map $h=q_{w_{0}} \circ q_{v_{0}}$, $\gamma_{1}$ has a superattracting periodic orbit of period $2n$. Let $U_{h}=U_{h}(\gamma_{1})$ be the immediate attracting basin of $\gamma_{1}$. Hence, if $K(h)$ is the filled Julia set of $h$, then $U_{h} \subset K(h)$ is a simply connected bounded open neighbourhood of $\gamma_{1}$ that is carried to itself by $h^{\circ n}$. We point out the two cases that could appear, depending on the behavior of the other two (complex) critical points of $h$, called $C_{1}$ and $C_{3}$.\\

\noindent \textbf{Case 1.} The map $h$ is hyperbolic (i.e. $C_{1}$ and $C_{3}$ are attracted).\\ 

\proof{ Each hyperbolic component in $P^{Q}$ is biholomorpfic to a Blaschke model. Within each of these components, the locus of the maps with a specific superattracting orbit is a smooth complex manifold. Each bones intersection is a center point for some hyperbolic component, and it has been proved that these intersections are transverse.}

\vspace{5mm}

\noindent \textbf{Case 2.} The map $h$ is not hyperbolic (i.e. $C_{1}$ and $C_{3}$ are not attracted to attracting cycles).\\

 \proof{ We will use quasiconformal surgery in the neighbourhood of our fixed map $h \in P^{Q}$. No iterates of the other two critical points of $h$ belong to $U_{h}$, the immediate attracting basin of $\gamma_{1}$, hence $U_{h}$ is isomorphic to the open unit disc, parametrized by its Bottcher coordinate. I.e., there exists a biholomorphic isomorphism that conjugates $h^{\circ n}$ to the squaring map $z \; \longrightarrow \; z^{2}$:\\

   $ \beta : U_{h} \; \longrightarrow \; \mathbb{D} $\\
 
   $ \beta(h^{\circ n}(z))=(\beta(z))^{2} $ \\

   We want to replace the superattracting basin $U_{h}$ by a basin with small positive multiplier $\Lambda$. For each $\Lambda$ in a small disc centered at zero, we will construct a new map $h_{\Lambda}$ corresponding to a $(v_{\Lambda},w_{\Lambda}) \in P^{Q}$ in such a way that $\Lambda \; \longrightarrow \; h_{\Lambda} \sim (v_{\Lambda},w_{\Lambda})$ is analytic and that $h_{0}=h$.

   The composition of smooth (analytic) maps

   $$\Lambda \; \longrightarrow \; h_{\Lambda} \sim (v_{\Lambda},w_{\Lambda}) \in P^{Q} \; \longrightarrow \; m(h_{\Lambda})$$

\noindent
is the identity. (Here $m$ denotes again the function that assigns to each map in $P^{Q}$ its multiplier at the specified attracting point). It follows that the partial derivatives $\frac{\partial m}{\partial v},\frac{\partial m}{\partial w}$ can't be simultaneously zero on a small neighbourhood of $h \in P^{Q}$. By the Implicit Function Theorem, the bone curve is smooth ${\cal{C}}^{1}$ on a small neighbourhood of $h$.}

\subsection{Quasiconformal surgery construction} \label{section:subsection2}
\noindent

    Consider the map $f(z)=z^{2}$ on the open unit disk $\mathbb{D}$ (which is the Bottcher parametrization of $h^{\circ n}$). Its unique critical point is the origin. Fix a small $\epsilon > 0$ (along the proof we will make specific requirements of how small we want $\epsilon$ to be) and let $\Lambda$ be an arbitrary complex number such that $0 \leq \mid \Lambda \mid \leq \epsilon$.

   Using a partition of unity, we perturb the map $f$ to a new degree 2 map $g_{\Lambda}$ such that:\\

   $\bullet$ $g_{\Lambda}$ has the same dynamics as $f_{\Lambda}(z)=z^{2}+\Lambda z$ inside a small disc around zero; in particular, the origin will be fixed, with multiplier $\Lambda$;\\

   $\bullet$ $g_{\Lambda}$ has the same dynamics as $f(z)=z^{2}$ outside a larger disc around zero.\\

   Choose a radius $r$ such that:

    $$ \frac{\epsilon}{2} \leq r \leq min(\frac{1}{2},1-\epsilon) $$
This will insure that $f_{\Lambda}$ maps $\Delta_{r^{2}}$ into itself and that the critical point of $f_{\Lambda}$ is in $\Delta_{r^{2}}$.

   Construct a ${\cal{C}}^{1}$ partition of unity $\rho:\mathbb{C} \; \longrightarrow \; \mathbb{R}$ with 

   $\bullet$ $\rho=0$ outside $\Delta_{\frac{r}{2}}$;
   
   $\bullet$ $\rho=1$ inside $\Delta_{r^{2}}$;

   $\bullet$ $0 \leq \rho \leq 1$ on $\Delta_{r/2} \backslash \Delta_{r^{2}}$\\

   Define $g_{\Lambda}:\mathbb{C} \; \longrightarrow \; \mathbb{C}$ as:

   $$g_{\Lambda}(z)=z^{2}+\Lambda \rho(z) z$$

   If we ask for $\frac{r}{2}(\frac{r}{2}+\epsilon) \leq r^{2}$, i.e.$\frac{2 \epsilon}{3} \leq r$ and by making $\epsilon$ smaller, if necessary, we can insure that $g_{\Lambda}$ has no critical point outside $\Delta_{r^{2}}$, for any $0 \leq \mid \Lambda \mid \leq \epsilon$. (Recall that the critical point of $g_{0}(z)=f(z)=z^{2}$ is $0 \in \Delta_{r^{2}}$ and the dependence $\Lambda \; \longrightarrow \; g_{\Lambda}$ is smooth for $\mid \Lambda \mid \leq \epsilon$). 

   \textbf{For short}: For any fixed $\mid \Lambda \mid \leq \epsilon$, the map $g_{\Lambda}:\mathbb{D} \; \rightarrow \; \mathbb{D}$ constructed above is a 2-to-1 ${\cal{C}}^{1}$ smooth map that carries $\Delta_{r} \backslash \Delta_{r^{2}}$ into $\Delta_{r^{2}}$ and carries $\Delta_{r^{2}}$ into itself. $g_{\Lambda}$ coincides with $f_{\Lambda}$ inside $\Delta_{r^{2}}$ and with $f$ outside of $\Delta_{\frac{r}{2}}$  (in particular it is conformal outside $\Delta_{r}$) and has no critical points in $\Delta_{r} \backslash \Delta_{r^{2}}$. We would like to emphasize that, as $\Delta_{r} \backslash \Delta_{r^{2}}$ is mapped by $g_{\Lambda}$ directly into $\Delta_{r^{2}}$, the annulus $\Delta_{r} \backslash \Delta_{r^{2}}$ is intersected at most once by any orbit under $g_{\Lambda}$.

   We pull $g_{\Lambda}$ back to $U_{h}$ through the Bottcher biholomorpic diffeomorphism $\beta$:

   $$ G_{\Lambda}=\beta^{-1} \circ g_{\Lambda} \circ \beta \, : U_{h} \; \rightarrow \; U_{h} $$

    The new map $G_{\Lambda}$ is 2-to-1 and ${\cal{C}}^{1}$ smooth, and has similar properties as the ones stated above for $g_{\Lambda}$ (see figure):

\begin{figure}[h]
\begin{center}
\input{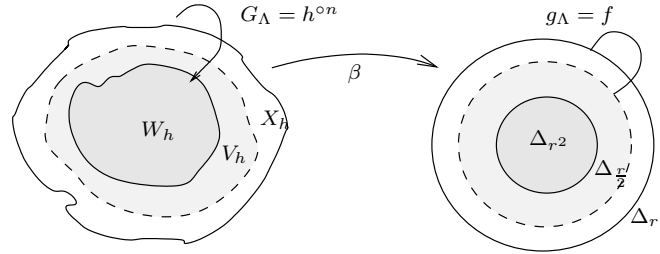}\\
\caption{\emph{$X_{h}$, $V_{h}$ and $W_{h}$ are the preimages under the Bottcher map $\beta$ of $\Delta_{r}$, $\Delta_{\frac{r}{2}}$ and $\Delta_{r^{2}}$, respectively. The map $g_{\Lambda}: \, \mathbb{D} \, \rightarrow \, \mathbb{D}$ pulls back as the ${\cal{C}}^{1}$-map $G_{\Lambda}$, that acts as $h^{\circ n}$ outside $V_{h}$ and carries $V_{h}$ to $W_{h}$.}}
\end{center}
\end{figure}

   But $h: \mathbb{C} \; \rightarrow \; \mathbb{C}$ carries

   $$U_{h} \, \rightarrow \, h(U_{h}) \, \xrightarrow{\sim} \, ...\, \xrightarrow{\sim} h^{\circ(n-1)}(U_{h}) \, \xrightarrow{\sim} \, h^{\circ n}(U_{h})=U_{h}$$
(acting as a diffeo except on $U_{h}$). So we can define $H_{\Lambda}$ as: 
  
   $H_{\Lambda}=h$ outside $V_{h}$ and\\ 

   $H_{\Lambda}=h^{^{\circ (1-n)}} \circ G_{\Lambda}$ inside $X_{h}$\\

   The new $H_{\Lambda}$ is ${\cal{C}}^{1}$ (notice that the two definitions coincide on $X_{h} \backslash V_{h}$) and has the desired dynamical behavior. However, it may fail to be analytic, hence it may not be a map in $P^{Q}$. The rest of the construction aims to transform $H_{\Lambda}$ into a polynomial $h_{\Lambda}\in P^{Q}$, preserving the dynamics.
   
   The Beltrami dilatation of $H_{\Lambda}$ is: 

   $$\mu_{H_{\Lambda}}(z)=\frac{(H_{\Lambda})_{\overline{z}}}{(H_{\Lambda})_{z}}$$
  
    Recall that $g_{\Lambda}$ has no critical point in $\Delta_{r} \backslash \Delta_{r^{2}}$, so $(g_{\Lambda})_{z} \neq 0$ on $\Delta_{r} \backslash \Delta_{r^{2}}$. Hence the denominator of:

   $$ \mu_{g_{\Lambda}}(z)=\frac{(g_{\Lambda})_{\overline{z}}}{(g_{\Lambda})_{z}}$$
never vanishes. Moreover, for fixed $z$, both top and bottom above are linear in $\Lambda$, so it follows easily that:

   $$ \Lambda \, \rightarrow \, \mu_{g_{\Lambda}}$$ is an analytic dependence. Hence $\mu_{H_{\Lambda}}(z)$ depends itself analytically on $\Lambda$ and :\\

   $\bullet$ $ \mu_{H_{\Lambda}}(z)=\mu_{G_{\Lambda}}(z)=\mu_{g_{\Lambda}}(\beta(z)) \frac{\overline{\beta'(z)}}{\beta(z)}$ on $X_{h}$ 

   $\bullet$ $\mu_{H_{\Lambda}}(z)=0$ outside $V_{h}$\\

   Under iteration of $g_{\Lambda}$, points hit the annulus $\Delta_{r} \backslash \Delta_{r^{2}}$ at most once, hence $\mu_{g_{\Lambda}}$ is bounded less than $1$ in modulus.

   $\mid \mu_{H_{\Lambda}}(z) \mid = \mid \mu_{g_{\Lambda}}(\beta(z)) \mid \mid \frac{\overline{\beta'(z)}}{\beta'(z)} \mid = \mid \mu_{g_{\Lambda}}(\beta(z)) \mid \leq 1$ on $X_{h} \backslash W_{h}$ and\\

   $\mu_{H_{\Lambda}}(z)=0$ outside $X_{h} \backslash W_{h}$.\\

   We define an ellipse field starting with circles inside $W_{h}$ and outside all preimages of $X_{h}$ under $H_{\Lambda}$ and pulling it back invariantly under $H_{\Lambda}$. All orbits hit $X_{h} \backslash W_{h}$ (the annular region where $H_{\Lambda}$ is not analytic) at most once, so the ellipse field is distorted at most once along any orbit. Let $\mu_{\Lambda}$ be the coefficient of this field. The dependence of $\mu_{\Lambda}$ on $\Lambda$ is holomorphic on $\mid \Lambda \mid \leq \epsilon$.

   Let $\phi_{\Lambda}$ solve the Beltrami equation:

   $$ \frac{\phi_{\overline{z}}}{\phi_{z}} = \mu_{\Lambda} $$ 
determined uniquely by the normalization $\phi_{\Lambda}(0)=0$, $\phi_{\Lambda}(1)=1$, $\phi_{\Lambda}(\infty)=\infty$,

   With this choice for $\phi_{\Lambda}$, $h_{\Lambda}=\phi_{\Lambda} \circ H_{\Lambda} \circ \phi_{\Lambda}^{-1}$ is a quartic complex polynomial. Moreover, for $\Lambda \in \mathbb{R}, \, \mid \Lambda \mid < \epsilon$, $h_{\Lambda}$ corresponds to a pair in the $Q$
-family (see ~\cite{R}). 
\subsection{The impossibility of bone-loops} \label{section:subsection3}
\indent

   Our plan for this section is to prove that bones in the parameter space $P^{Q}$ can not contain any loops (i.e. simple closed curves). Recall that we proved in section 2 that each bone contains a simple bone-arc connecting two boundary points, and that all possible distinguished kneading data of the bone can be found in a certain order along this bone-arc.

   We argue by contradiction. Suppose there exists a bone loop $L$. We will show next that the interior ${\cal{U}}$ of the loop can't contain any hyperbolic maps. This will contradict the genericity of hyperbolicity stated in theorem 3.7.\\

   \noindent \textbf{Remark.} The following statements and proofs are given for left bones, but apply by symmetry to right bones.\\

  \begin{lemma} A left bone loop in $P^{Q}$ can't contain any distinguished point, hence it can't contain any crossing with a right bone.
  \end{lemma}

   \proof{ Any distinguished point on the loop $L$ would need to have a kneading-data already achieved along the bone arc. Thurston's Theorem shows easily that this is impossible.}

   \begin{thm} The region enclosed by a left bone loop in $P^{Q}$ can't contain any hyperbolic maps.
   \end{thm}

   \proof{ We know by theorem 3.6 that each hyperbolic component in $P^{Q}$ is an open topological 2-cell that contains a unique post-critically finite point, called \emph{``center''}. Moreover, the intersection of any bone with a hyperbolic component must be a simple arc passing through the center.

   Suppose, by contradiction, that some hyperbolic component ${\cal{H}}$ intersects the region ${\cal{U}}$. We have two cases:\\

   (1) ${\cal{H}} \subset {\cal{U}}$. then there is a bone that passes through the center of ${\cal{H}}$. This can only be a bone arc, as bone loops can't contain distinguished points (by lemma 4.2). From the Jordan Curve Theorem, this bone arc has to intersect the bone loop $L$, contradiction with lemma 4.2.\\

   (2) ${\cal{H}}$ intersects the loop $L$. Then the loop must contain the center point of ${\cal{H}}$, again contradiction. }

\section{Topological conclusions}
\subsection{The entropy and the bones} \label{section:subsection1}
\vspace{5mm}
\indent

   Recall that our final claim is: for each fixed $h_{0} \in [0,\log4]$, the level-set $i(h_{0})=\{ h=h_{0} \}$ of the entropy function in either parameter space, called $h_{0}$-isentrope, is connected.

   In the $ST$-family, the analysis of the properties of entropy level-sets is an easy exercise. One can obtain the following fairly straight-foreward (see ~\cite{MT} and ~\cite{R}):

   \begin{thm} In $P^{ST}$, the entropy is a monotone function of either coordinate. For each $h_{0} \in [0,\log4]$, the corresponding $h_{0}$-isentrope is contractible, as it is a deformation retract of the contractible region $\{h \leq h_{0}\}$.
   \end{thm}

   To obtain similar results in the quartic family, we will need some notations and results from the general theory of $m$-modal maps of the interval.

   If $f: I \rightarrow I$ is an $m$-modal map with folding points $c_{1} \leq c_{2} \leq ... \leq c_{m}$, we define the sign of the fixed point $x$ of $f^{\circ k}$ with itinerary $\Im(x)=(A_{0},A_{1},...,A_{k-1})$ as the number:

   $$ sign(x) = \epsilon(A_{0}) \epsilon(A_{1}) ... \epsilon(A_{k-1}) $$
where $\epsilon(A_{j})=+1, \; -1$ or $0$ according to $A_{j}$ being an increasing/decreasing lap of $f$ or a folding point $c_{1},...,c_{m}$. If $sign(x)=-1$ we say that $x$ is a fixed point of negative type of $f^{\circ k}$.

   We define $Neg(f^{\circ k})$ as the number of fixed points of negative type of $f^{\circ k}$.

  \begin{thm} (\emph{[MT], page 22}) If $f$ is an interval $m$-modal map, then its topological entropy is:

  \[h(f) = \overline{\lim_{k\rightarrow \infty}} \frac{1}{k} \log^{+}(Neg (f^{\circ k}))\]
where $\log^{+}s = max( \log(s),0 )$.
  \end{thm}

  \noindent \textbf{Remark:} $Neg (f^{\circ k})$ is an integer $\geq 1$ unless $f^{\circ k}$ has no fixed points of negative type; in that case, $ \log^{+}(Neg (f^{\circ k})) = 0 $.\\

   The following result is a simple consequence of theorem 5.2 (see ~\cite{R} for proof and details).

   \begin{lemma} If for two $m$-modal interval maps $f$ and $g$ the topological entropies $h(f) \neq h(g)$, then the sequence $\mid Neg (f^{\circ k}) - Neg (g^{\circ k}) \mid$ must be unbounded as $k \rightarrow \infty$.
   \end{lemma}

\noindent \textbf{Notation.} For $p=(v,w) \in P^{Q}$, call $Q_{p}=q_{w} \circ q_{v}$ and for $p=(v,w) \in P^{ST}$, call $ST_{p}=st_{w} \circ st_{v}$.\\

   \begin{lemma} Consider $p_{1}=(v_{1},w_{1})$ and $p_{2}=(v_{2},w_{2})$ in $P^{Q}$ such that
   
   $$ h(Q_{p_{1}}) \neq h(Q_{p_{2}}) $$
   Then any path in $P^{Q}$ from $p_{1}$ to $p_{2}$ crosses infinitely many bones.
   \end{lemma}

   \proof{ 

   Consider an arbitrary path in $P^{Q}$ from $p_{1}$ to $p_{2}$:\\

   $\quad \quad  p:[0,1] \rightarrow P^{Q} \; , \; p(t) = (v(t),w(t)) $

   $\quad \quad p(0) = p_{1} = (v_{1},w_{1}) \; , \; p(1) = p_{2} = (v_{2},w_{2}) $ \\

   For a fixed $k \in \mathbb{N}$, as $t$ goes from $0$ to $1$, $Neg (Q_{p(t)}^{\circ k})$ changes whenever a fixed point of $Q_{p(t)}^{\circ k}$ ( i.e. a periodic point of $Q_{p(t)}$ of period dividing $k$ ) of negative type appears or disappears. An existing negative-type fixed point of $Q_{p(t)}^{\circ k}$ can be lost under continuous deformations of the map by becoming a positive-type fixed point. Conversely, a such fixed point can appear by a reverse process. Both changes imply the existence of an intermediate state, corresponding to some $t^{*} \in [0,1]$, in which the respective fixed point is a critical point of $Q_{p(t^{*})}^{\circ k}$. 

   In other words, a critical point of $Q_{p(t^{*})}$ has to be periodic of period dividing $k$. This implies that $p(t^{*}) = (v(t^{*}),w(t^{*})) \in P^{Q}$ is on either a left or a right bone of period $2n \mid 2k$.

   So if the integer $Neg \, ((q_{w(t)} \circ q_{v(t)})^{\circ k})$ has an actual change at $t=t^{*}$, then the path $p(t)$ crosses a bone at $t=t^{*}$.

   To end the proof of the lemma, suppose that the path $p(t)$ only crosses $N$ bones. Then, for all $k \in \mathbb{N}$,

   $$ \mid Neg (Q_{p_{1}}^{\circ k}) - Neg (Q_{p_{2}}^{\circ k}) \mid$$
would be bounded by $N$, contradiction with lemma 5.3.  }

\subsection{The entropy and the cellular structure} \label{section:subsection2}
\vspace{5mm}
\indent

  Recall that either parameter space $P^{ST}$ and $P^{Q}$ has for each fixed value of $n$ an associated cellular complex structure, called $P^{ST}_{n}$ and $P^{Q}_{n}$, respectivelly. The two cell complexes are homeomorphic through the function $\eta$ defined in section 2.11.
   
   The following lemma is valid for either complexes $P_{n}=P^{ST}_{n}$ or $P_{n}=P^{Q}_{n}$.

   \begin{lemma} For any $\epsilon > 0$, there exists $n \in \mathbb{N}$ such that, if $p$ and $p'$ belong to the same closed cell in $P^{n}$, then the corresponding maps satisfy:

   $$ \mid h_{p} - h_{p'} \mid < \epsilon $$
   \end{lemma}

   \proof{ Suppose the contrary: there exists $\epsilon > 0$ such that, for all $n \in \mathbb{N}$, there are two parameters $p_{n}$ and $p_{n}'$ in some common cell of $P_{n}$ with:

   $$ \mid h_{p_{n}} - h_{p'_{n}} \mid \geq \epsilon $$

   By the compactness of $P$, we can choose a subsequence $(k_{n})_{n} \subset \mathbb{N}$ such that both $(p_{k_{n}})_{n}$ and $(p'_{k_{n}})_{n}$ converge in $P$:

   $ \hspace{20mm} p_{k_{n}} \longrightarrow p \quad$ as $n \rightarrow \infty$

   $ \hspace{20mm} p'_{k_{n}} \longrightarrow p' \quad$ as $n \rightarrow \infty$

   The entropy function is a continuous function of parameters in either family (see for example ~\cite{MT}). Using this and passing to the limit:

   $$ \mid h_{p_{k_{n}}} - h_{p'_{k_{n}}} \mid \geq \epsilon \quad \Rightarrow \quad \mid h_{p} - h_{p'} \mid \geq \epsilon $$

   Moreover, the closed cells of $P_{n}$ are nested as $n$ increases (in other words, the cell complex gets \emph{``finer''} with larger values of $n$ ).

   Fix an arbitrary $N \in \mathbb{N}$. For all $k_{n} \geq N$, $p_{k_{n}}$ and $p'_{k_{n}}$ are in the same closed cell of $P_{k_{n}}$, hence in the same closed cell of $P_{N}$. 

   In conclusion, for any arbitrary $N \in \mathbb{N}$, $p$ and $p'$ are in the same closed cell of $P_{N}$, yet:

   $$ \mid h_{p} - h_{p'} \mid \geq \epsilon > 0 $$
contradiction with lemma ... }

   \begin{lemma} Fix $n \in \mathbb{N}$. In either parameter space $P$, the entropy function:

    $$ P_{n} \, \longrightarrow [0,log4] $$
    $$ p \, \longrightarrow h(g_{p}) $$
restricted to any closed cell in $P_{n}$ takes its maximum and minimum values on the boundary of the cell ( more precisely on the boundary vertexes ).
   \end{lemma}

   \proof{ In the case $P_{n}=P^{ST}_{n}$, the proof is a simple corollary of lemma ... We have to prove the identical statement for $P_{n}=P^{Q}_{n}$.

   For the fixed $n \in \mathbb{N}$, suppose the lemma is not true for some closed cell $C^{Q}_{n} \in P^{Q}_{n}$, that is : there exists $p^{*}=(v^{*},w^{*}) \in int(C^{Q}_{n})$ such that 

  $$ h(Q_{p^{*}}) = h(q_{w^{*}} \circ q_{v^{*}}) > h_{max} \: , $$
where $h_{max}$ is the maximum value of the entropy on the boundary $\delta(C^{Q}_{n})$.

   Let $$\epsilon = \frac{h(Q_{p^{*}}) - h_{max}}{2} \geq 0 $$

   By lemma 5.5, there exists $m \in \mathbb{N}$ such that the entropy variation on all closed cells of $P^{Q}_{m}$ is less than $\epsilon$. WLOG, we can take $m>n$. Call $C^{Q}_{m}$ the closed cell in $P^{Q}_{m}$ such that $p^{*} \in C^{Q}_{m} \subset C^{Q}_{n}$ and consider any arbitrary vertex $p_{m} = (v_{m},w_{m})$ of $C^{Q}_{m}$.

   As $p^{*}, \: p_{m} \in C^{q}_{m}$, we automatically have:

   $$ \mid h(Q_{p^{*}}) - h(Q_{p_{m}}) \mid < \epsilon $$

   But $ \quad h_{max} + 2\epsilon \leq h(Q_{p^{*}}) \quad $, so:

   $$ h(Q_{p_{m}}) > h_{max} $$

   The homeomorphism of complexes $\quad \eta_{m}^{-1}:P^{Q}_{m} \longrightarrow P^{ST}_{m} \quad $ carries vertexes to vertexes with the same entropy, edge to edge with the same interval of entropies and 2-cells to 2-cells. So $ \quad C^{ST}_{m} = \eta_{m}^{-1}(C^{Q}_{m}) \quad $ will be a 2-cell in $P^{ST}_{m}$ and $ \quad q_{m} = \eta_{m}^{-1}(p_{m}) $ will be a vertex of $C^{ST}_{m}$. Also, $\quad \eta_{n}^{-1}(\delta C^{Q}_{n}) = \delta(C^{Q}_{n}) = \delta C^{ST}_{n} \quad$, so the maximum value $h_{max}(\delta C^{Q}_{n})$ of the entropy on $\delta C^{Q}_{n}$ is the same as the maximum value $h_{max}(\delta C^{ST}_{n})$ on $C^{ST}_{n}$. Hence, in the stunted family:

   $$ h(ST_{q_{m}}) = h(Q_{p_{m}}) > h_{max}(\delta C^{Q}_{n}) = h_{max}(\delta C^{ST}_{n}) \:,$$
contradiction, since the result has already been proved for $P^{ST}$. }

  \begin{corol} For a fixed $n \in \mathbb{N}$, the interval of entropy values realized by any cell in $P^{Q}_{n}$ is the same as the interval of values for the corresponding cell in $P^{ST}_{N}$.
  \end{corol}

   For $h_{0} \in [0,log4]$ we will use the following notation for the $h_{0}$-isentrope in either family:\\

   $ i^{ST}(h_{0}) = \{ (v,w) \in P^{ST} \; \slash \; h(st_{w} \circ st_{v}) = h_{0} \}$\\

   $  i^{Q}(h_{0}) = \{ (v,w) \in P^{Q} \; \slash \; h(q_{w} \circ q_{v}) = h_{0} \}$\\

   For a fixed $n \in \mathbb{N}^{*}$, we also use the following notations:\\

   $ N^{ST}_{n}(h_{0}) = \bigcup { \{ C^{ST}_{n} \; \slash \; C^{ST}_{n} \in P^{ST}_{n}, \; C^{ST}_{n} \cap i^{ST}(h_{0}) \neq \Phi \} } $\\

   $ N^{Q}_{n}(h_{0}) = \bigcup { \{ C^{Q}_{n} \; \slash \; C^{Q}_{n} \in P^{Q}_{n}, \; C^{Q}_{n} \cap i^{Q}(h_{0}) \neq \Phi \} } $\\

   \noindent \textbf{Remarks :} (1) Clearly: $ i^{ST}(h_{0}) \subset N^{ST}_{n}(h_{0}) $ and $ i^{Q}(h_{0}) \subset N^{Q}_{n}(h_{0}) $.

   (2) Recall that for fixed $n$ we have the homeomorphism of cell complexes:
       
    $$ \eta_{n} : P^{ST}_{n} \rightarrow  P^{Q}_{n} $$

    If $C^{ST}_{n}$ is a cell in $P^{ST}_{n}$ that touches $i^{ST}(h_{0})$, then the corresponding cell $C^{Q}_{n} = \eta_{n}(C^{ST}_{n})$ will touch $i^{Q}(h_{0})$ and conversely. This follows from corollary 5.7, which states that the interval of entropy values is the same in the two closed cells $C^{ST}_{n}$ and $C^{Q}_{n}$.

    Fix an entropy value $h_{0} \in [0,log4]$ and an $n \in \mathbb{N}^{*}$.

    Since $N^{ST}_{n}(h_{0})$ and $N^{Q}_{n}(h_{0})$ are both unions of closed cells, they are compact subsets of $P^{ST}$ and $P^{Q}$, respectively. By the previous theorem, $N^{ST}_{n}(h_{0})$ is connected, so its image $N^{Q}_{n}(h_{0}) = \eta_{n}(N^{ST}_{n}(h_{0})$ is also connected. Hence we have the following:\\

\begin{figure}[h]
\begin{center}
\input{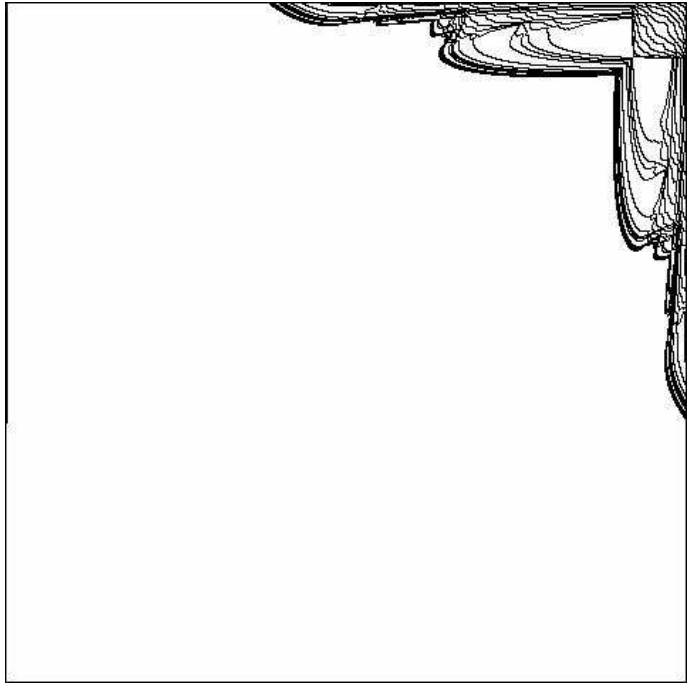}\\
\caption{\emph{ The isentropes in $P^{Q}$ appear to be either arcs joining two points in $\partial P^{Q}$, or connected regions between such arcs, or a single point (the case $(v,w)=(1,1)$ of entropy $\log 4$.}}
\end{center}
\end{figure}

   \noindent \textbf{Summary.} \emph{For any $n \in \mathbb{N}^{*}$, the set $N^{Q}_{n}(h_{0})$ is compact, connected and contains $i^{Q}(h_{0})$}.\\

    We have now a quite comprehensive description of the sets $N^{Q}_{n}(h_{0})$. To obtain topological properties of $i^{Q}(h_{0})$, we try to relate it to the collection $\{ N^{Q}_{n}(h_{0}) \}_{n \in \mathbb{N}}$.

   \begin{lemma} $\bigcap {N^{Q}_{n}(h_{0})} = i^{Q}(h_{0}) $
   \end{lemma}

   \proof{ Since $i^{Q}(h_{0}) \subset N^{Q}_{n}(h_{0})$ for all $n \in \mathbb{N}^{*}$, the inclusion $i^{Q}(h_{0}) \subset \bigcap N^{Q}_{n}(h_{0})$ is trivial.

   For the converse, suppose there exists $(v,w) \in \bigcap N^{Q}_{n}(h_{0}) \backslash i^{Q}(h_{0})$. In other words: for any arbitrary $n \in \mathbb{N}^{*}$, $(v,w)$ is contained in a closed cell $C^{Q}_{n} \subset P^{Q}_{n}$ that touches $i^{Q}(h_{0})$, but such that $(v,w) \notin i^{Q}(h_{0})$. For any such closed cell $C^{Q}_{n}$, there exists $(v_{n}^{*},w_{n}^{*}) \in i^{Q}(h_{0}) \cap C^{Q}_{n}$.

   The sequence $(v_{n}^{*},w_{n}^{*})_{n \in \mathbb{N}^{*}}$ satisfies in particular:\\

   (1) $ (v_{n}^{*},w_{n}^{*}) \neq (v,w), \; \forall n \in \mathbb{N}^{*} $ \\
 
   (2) $ h(q_{w_{n}^{*}} \circ q_{v_{n}^{*}}) = h_{0} $\\

   We calculate:

   $$ \mid h(q_{w^{*}} \circ q_{v^{*}}) - h(q_{w} \circ q_{v}) \mid \, = \, \mid h_{0} - h(q_{w} \circ q_{v}) \mid $$

   This contradicts the statement of lemma 5.5: the maximal variation of the entropy over cells in $P^{Q}_{n}$ can be made arbitrarily small by increasing $n$. }

   \begin{thm} $i^{Q}(h_{0})$, the $h_{0}$-isentrope in $P^{Q}$, is connected.
   \end{thm}

   \proof{ $i^{Q}(h_{0})$ is an intersection of compact, connected sets in $P^{Q}$, therefore it is compact and connected.}

\bibliographystyle{plain}

\end{document}